\documentclass[a4paper]{article} %version ecran

\usepackage[english]{babel}

\usepackage[utf8]{inputenc}
\usepackage{amsmath,amsthm}
\usepackage{amsfonts,amssymb}

\usepackage{hyperref}
\usepackage{graphicx}

\usepackage{tikz,tikz-cd} %Graĥique

\usepackage{xifthen} %provide "\ifthenelse" et "\isempty{}"

\usepackage{titlesec}
\titleformat{\subsubsection}[runin]{\normalfont}{\thesubsubsection}{0pt}{}[.]

\newcommand{\toposdefaut}{0}
\newcommand{\topos}[1][\toposdefaut]{ 
\ifthenelse{\equal{#1}{0}}{ \mathcal{T} }
{
\ifthenelse{\equal{#1}{1}}{ \mathcal{E} }{ #1 }
}
}

\newcommand{\sh}{\textsf{Sh}}

\newcommand{\fn}{\text{Fn }}

\newcommand{\spec}{\text{Spec }}

\newcommand{\Q}{\mathbb{Q}}
\newcommand{\R}{\mathbb{R}}
\newcommand{\C}{\mathbb{C}}
\newcommand{\N}{\mathbb{N}}
\newcommand{\Z}{\mathbb{Z}}

\newcommand{\Ecal}{\mathcal{E}} 
 
\newcommand{\Tcal}{\mathcal{T}}

\newcommand{\Ocal}{\mathcal{O}}

\newcommand{\Scal}{\mathcal{S}} 
 
\newcommand{\Fcal}{\mathcal{F}} 
 
\newcommand{\Hcal}{\mathcal{H}}

\newcommand{\Lcal}{\mathcal{L}} 
\newcommand{\Mcal}{\mathcal{M}}

\newcommand{\Ccal}{\mathcal{C}} 
 
\newcommand{\Bcal}{\mathcal{B}}

\newcommand{\Length}{\overleftarrow{\mathbb{R}_+^\infty}}

\newcommand{\block}[1]
{
\par \subsubsection{} #1

\bigskip}

\newcommand{\blockn}[1]{\par #1 \bigskip}

\newcommand{\Th}[1]
	{
	\bigskip	
	\textbf{Theorem : }{\itshape #1}
		
	\bigskip
	}

\newcommand{\Prop}[1]
	{

	\bigskip
	
	\textbf{Proposition : }{\itshape #1}
		
	\bigskip
	
	}

\newcommand{\Cor}[1]
	{

	\bigskip
	
	\textbf{Corollary : }{\itshape #1}	
		
	\bigskip

	}

\newcommand{\Lem}[1]
	{

	\bigskip
	
	\textbf{Lemma : }{\itshape #1}
		
	\bigskip
	
	}

\newcommand{\Def}[1]
	{
	
	\bigskip
	
	\textbf{Definition : }{\itshape #1}
	
	\bigskip
	
	}

\newcommand{\Dem}[1]{
	
	\smallskip
	
	\textbf{Proof : } \par
	 {#1} $\square$

	 \bigskip
}

\begin{document}

\selectlanguage{english}
\pagestyle{plain}
\title{Localic metric spaces and the localic Gelfand duality}
\author{Simon Henry}
%\address{Simon.~Henry: \\ Paris, F-75005 France}

\renewcommand{\thefootnote}{\fnsymbol{footnote}} 
\footnotetext{\emph{Keywords.} Locales, Metric locales, Banach locales, Gelfand duality.}
\footnotetext{\emph{2010 Mathematics Subject Classification.} 18B25, 03G30, 06D22, 46L05, 47S30.}
\footnotetext{\emph{email:} shenry2@uottawa.ca}
\renewcommand{\thefootnote}{\arabic{footnote}} 

% 18B25 : topoi
% 03G30 : Categorical logic, topoi
% 06D22 : Frames, Locales
% 46L05 : General theory of $C^*$-algebras
% 47S30 : Constructive operator theory 

\maketitle

\begin{abstract}

In this paper we prove, as conjectured by B.Banachewski and C.J.Mulvey, that the constructive Gelfand duality can be extended into a duality between compact regular locales and unital abelian localic $C^{*}$-algebras. In order to do so we develop a constructive theory of localic metric spaces and localic Banach spaces, we study the notion of localic completion of such objects and the behaviour of these constructions with respect to pull-back along geometric morphisms.
\end{abstract}
  
\tableofcontents

\section{Introduction}

\blockn{In \cite{banaschewski2006globalisation}, C.J.Mulvey and B.Banaschewski showed\footnote{To be more accurate, they only showed this result internally in Grothendieck toposes, using at some points an external argument relying on the axiom of choice (the Barr covering theorem). A completely internal and constructive proof has been given later by T.Coquand and B.Spitters in \cite{coquand2009constructive}.} that the usual Gelfand duality between abelian $C^*$ algebras and compact (Hausdorff) topological spaces can be extended into a ``constructive" Gelfand duality between $C^*$ algebras and compact completely regular locales. A locale (see \ref{introlocale}) is almost the same as a topological space, but may fail to have points. A locale which has enough points is called a spatial locale and is the same thing as a (sober) topological space. Assuming the axiom of choice, any locally compact locale has enough points; hence the result of Banaschewski and Mulvey gives back the usual Gelfand duality when assuming the axiom of choice. But the constructive version can be applied to a broader context: for example an internal application to the topos of sheaves over a topological space $X$ relates continuous fields of abelian $C^*$ algebras over $X$ and proper maps to $X$, and this can also be applied to more general toposes.

}

\blockn{At the end of their proof of the constructive Gelfand duality, Banachewski and Mulvey suggested that ``compact completely regular" is not the most natural condition one would have expected. It would be nicer to weaken this condition into ``compact regular" (which is the same as compact separated, see \cite{sketches} C3.2.10). Unfortunately, when a locale is not completely regular it might fail to have $\C$-valued continuous functions, and hence the associated $C^*$ algebra has no reason to keep track of enough informations about $X$. They suggested that their result should be extended into a duality between compact regular locales and a notion of localic $C^*$ algebras yet to be defined. This is a natural idea because when $X$ is a compact regular locale, one can still define a locale $[X,\C]$ of functions from $X$ to $\C$ and complete regularity only concerns the existence of points for this locale.
The main goal of this article is to define this notion of localic $C^*$ algebras (which we will call $C^*$ locales) and to prove this conjectured duality.}

\blockn{Two other reasons for developing a theory of localic $C^*$ algebras and more generally of localic Banach spaces (called Banach locales) are the following. In \cite{moerdijk1990classifying} I.Moerdijk showed (using the result of A.Joyal and M.Tierney in \cite{joyal1984extension}) that Grothendieck toposes can be identified with a full subcategory of the $2$-category of localic groupoids (that is groupoids in the category of locale, morphisms between them being the localic principal bi-bundles, see \cite{bunge1990descent} for more details). A Banach space in the logic of the topos which corresponds to a localic groupoid $G_1 \rightrightarrows G_0$ is essentially a continuous field of Banach spaces $\Bcal$ over $G_0$ endowed with a continuous action of $G_1$ such that there are enough local sections of $\Bcal$ which have an open stabilizer. This hypothesis of open stabilizers is, from the point of view of analysis and geometry, a little too restrictive and is related to the requirement of existence of points. Hence one could expect that a good notion of Banach locale could remove it. Also for the purpose of non-commutative geometry one would like to be able to study equivariant bundle on general localic (topological) groupoids and not just those which correspond to toposes. For example the groupoid defined by $G_0$ being a point and $G_1$ being a connected locally compact topological group does not correspond to a topos but is an important groupoid for non-commutative geometry. In order to define the notion of Banach space over an arbitrary localic groupoid, an important point is that this notion should descend along open surjections (see \ref{Descenttheory}). Unfortunately, there is no such descent property for Banach spaces and $C^*$ algebras. However, as locales descend along open (or proper) surjections and as the pull-back of Banach spaces is a pull-back of the localic completion, we will be able to prove that Banach locales and $C^*$ locales have this descent property, and form in fact the ``stackification" of the notion of Banach spaces and $C^*$ algebras, i.e. the smallest generalization of the notion which have this descent property.}

\blockn{Section \ref{sectionnotationprelim} reviews some well known facts and definitions, mostly about the theory of locales, in order to fix the notation and prove some basic but important results for the rest of the paper. In section \ref{sectionmetriclocale} we will develop the theory of metric locales in a constructive context (the classical theory is already known and can be found for example in \cite{picado2012frames}). We also show how to construct a classifying locale $[X,Y]_1$ for metric maps between two complete metric locales, which was apparently not known even in the classical case. In section \ref{sectionBanach} we apply the theory of section \ref{sectionmetriclocale} in order to define Banach locales and $C^*$ locales and prove the announced result, although most of the technical difficulties lie in section \ref{sectionmetriclocale}.}

\blockn{An extended version of this article can be found in the author's thesis \cite{mythesis}. This extended version also contains an additional section were we prove (assuming the axiom of dependant choice) that when we work internally in a topos $\Tcal$ satisfying some technical condition related to paracompactness then the category of localic banach spaces of $\Tcal$ is equivalent to the category of usual Banach spaces of $\Tcal$. This result is a topos theoretic adaptation of a theorem\footnote{published as an appendix of \cite{fell1977induced}} of A.Douady and L.Dal Soglio-Herault asserting that over a paracompact topological spaces every Banach bundle has enough continuous sections. We decided not to include this last result in the present paper because we think that it still needs to be improved, in particular, more recent results we obtained suggest that it should be a consequence of a fully constructive result with more natural hypothesis.}

\section{Notations and Preliminaries}
\label{sectionnotationprelim}

\subsection{General remarks}

\blockn{In all the article we are implicitly working internally in an elementary topos $\Scal$ with a natural number object $\N$. This means that we will never use neither the law of excluded middle nor the axiom of choice. Objects of $\Scal$ will simply be called ``sets". All other toposes mentioned are bounded $\Scal$-toposes, i.e. Grothendieck toposes over $\Scal$ (although the hypothesis bounded can probably be removed most of the time).}

\blockn{A proposition (internal to a topos) is said to be \emph{decidable} if it is complemented (i.e. such that $P \vee \neg P$ holds). An object is said to have decidable equality, or to be \emph{decidable}, if its diagonal embedding $X \rightarrow X \times X$ is complemented.
A set $X$ (or an object of a topos) is said to be \emph{inhabited} if it satisfies (internally) $\exists x \in X$. It is said to be \emph{finite} if it is Kuratowski finite (see \cite[D5.4]{sketches}), i.e. if $\exists n \in \N, x_1, \dots x_n \in X$ such that $\forall x \in X, \exists i , x = x_i$. Note that in particular (as $\N$ is decidable) a finite set is either empty or inhabited.
}

\blockn{When considering product $E_1 \times \dots \times E_n$ of objects of any kind (generally locales) we will denote by $\pi_i$ the projection onto $E_i$, by $\pi_{i,j}$ the projection onto $E_i \times E_j$, etc... We generally do not specify the domain of definition and we hope that it will be clear from the context. For example one has: $\pi_1 \circ \pi_{i,j} = \pi_i$ and $\pi_2 \circ \pi_{i,j} = \pi_j$ because in these formulas $\pi_1$ and $\pi_2$ denote the two projections from $E_i \times E_j$ to $E_i$ and $E_j$ respectively.}

\subsection{The category of locales}
\label{introlocale}
\blockn{We will start by briefly introducing the notion of locale, essentially in order to fix the notation and the vocabulary. A short introduction to this subject can be found in the first two sections of \cite{borceux3}, a more complete one in part $C$ (especially in $C1$) of \cite{sketches} and an extremely complete (but non constructive) one in \cite{picado2012frames}.}

\block{A \emph{frame} is an ordered set which admit arbitrary supremums and such that binary infimums distribute over arbitrary supremums. A morphism of frame is a non-decreasing map which preserve both arbitrary supremum and finite infimum. }

\block{The category of \emph{locales} is defined as the opposite category of the category of frames. But we will adopt ``topological" notations for them:
\begin{itemize}
\item If $X$ is a locale, the corresponding frame is denoted by $\Ocal(X)$.
\item If $f : X \rightarrow Y$ is a morphism of locales, we denote by $f^*$ the corresponding frame homomorphism from $\Ocal(Y)$ to $\Ocal(X)$.
\item An element $U \in \Ocal(X)$ is called an open sublocale of $X$, the top element of $\Ocal(X)$ is denoted $X$.
\item As $f^*$ commutes to arbitrary supremums, it has a right adjoint denoted $f_*$.
\end{itemize}

Also we will tend to call unions and intersections the supremums and infimums in $\Ocal(X)$.

}

\block{A \emph{sublocale} of a locale $X$ is (an equivalence class of) a locale $Y$ endowed with a morphism $f :Y \rightarrow X$ such that $f^*$ is a surjective frame homomorphism (such a morphism is called an \emph{inclusion}). A morphism of locale $f$ is said to be \emph{surjective} if the corresponding frame homomorphism is injective. In particular, the injection/surjection factorisation of frame homomorphisms induces a unique (up to unique isomorphism) factorisation of every morphism of locale $f : X \rightarrow Y$ in a surjection followed by an inclusion:

\[ X \twoheadrightarrow f_!(X) \hookrightarrow Y. \]

The sublocale $f_!(X)$ is called the image\footnote{From a purely categorical point of view, we should call it the regular image of $X$.} of $f$. More generally if $S$ is any sublocale of $X$ we denote by $f_!(S)$ the image of the restriction of $f$ to $S$ and this is called the image of $S$ by $f$.
}

\block{If $f :X \rightarrow Y$ is a morphism of locales and $S$ is a sublocale of $Y$ then the categorical pull-back $f^{-1}(S)$ is a sublocale of $X$ and one has an adjunction formula:

\[ A \subset f^{-1}(B) \Leftrightarrow f_!(A) \subset B \]

for any sublocale $A$ of $X$ and $B$ of $Y$.
}

\block{If $U$ is an element of the frame $\Ocal(X)$ then it corresponds to a sublocale (also denoted $U$) of $X$ which is defined by the frame $\Ocal(U) = \{ v \in \Ocal(X) | v \leqslant U \}$ and which is sent into $X$ by the morphism corresponding to $i^*(V)=V \wedge U$ for any $V \in \Ocal(X)$. Hence, the elements of $\Ocal(X)$ correspond to particular sublocales of $X$, which justifies the term ``open sublocales" for elements of $\Ocal(X)$. Also, through this identification, one has $f^*(U)= f^{-1}(U)$. }

\block{To any locale $X$ one can associate the topos of sheaves on $X$, denoted $\sh(X)$. If $X$ and $Y$ are two locales, the category of geometric morphisms from $\sh(X)$ to $\sh(Y)$ is (equivalent to) the ordered set of locale morphisms from $X$ to $Y$ ordered by the pointwise ordering of the corresponding frame homomorphism (this is called the specialisation order). For this reasons locales will be seen as a specific kind of toposes.}

\block{An extremely important result of the theory of locales, that we will use constantly, is that there is an equivalence of category between $X$-locales, that is locales in the logic of $\sh(X)$ and locales $Y$ endowed with a morphism to $X$. This allows one to turn any reasonable property of locales into a property of geometric morphisms, corresponding to the relative notion, for example one says that a map $Y \rightarrow X$ is proper if the $X$-locale corresponding to $Y$ is compact in the logic of $\sh(X)$.
}

\block{At several points of this article we will deal (in simple situations) with locales as if they had points in order to define a map between two locales or to give constraints on some map. This kind of expression should of course not be interpreted in terms of points of a locale $X$ but in terms of ``generalized points", that is morphisms from $T$ to $X$ for an arbitrary locale $T$, and all the constructions done on these points should be interpreted in the logic of $\sh(T)$. If all the constructions on these generalized elements are ``geometric" (that is compatible with the pull-back from $\sh(T)$ to $\sh(T')$ for any locale $T'$ over $T$) then these constructions yield a morphism of functor, or relation between such morphisms and hence by the Yoneda lemma this indeed gives a morphism of locales or conditions between such morphisms. } 

\block{One says that a locale $\Lcal$ \emph{classifies} some theory $T$ if the topos $\sh(\Lcal)$ classifies the theory $T$. (see the part $D$ of \cite{sketches} for the general theory of classifying toposes) Locales are the classifying spaces of what is called propositional geometric theory. That is geometric theory over a signature (see \cite{sketches}D1.1.1) which contains no sorts. In particular it contains no function symbol and all the relations symbol it contains have no free variable and are called propositions. A propositional geometric theory classified by a locale $\Lcal$ is essentially the same thing as a presentation of the frame $\Ocal(\Lcal)$ : indeed basic proposition of the theory are generator and the geometric (in the sense of \cite[D1.1.3(f)]{sketches}) axioms of the theory are relations of the form $T \leqslant T'$ where $T$ and $T'$ are formed from the basic proposition using finite intersection and arbitrary union.
}

\block{\label{geompblocale}If $\Lcal$ is a locale in the logic of some topos $\Tcal$ and if $f : \Ecal \rightarrow \Tcal$ is a geometric morphism then, $f^* \Ocal(\Lcal)$ is in general not a frame in $\Ecal$, but it can be completed in a frame, giving rise to locale called $f^\#( \Lcal)$ in $\Ecal$. More precisely, if one takes a presentation of $\Lcal$, then one can pull-back the presentation through $f$ and construct a locale $\Lcal'$ in $\Ecal$. One can then check from the universal property that one has the following pull-back diagram of toposes:

\[
\begin{tikzcd}[ampersand replacement=\&]
\sh_{\Ecal}( \Lcal') \arrow{r} \arrow{d} \& \Ecal \arrow{d} \\
\sh_{\Tcal}(\Lcal) \arrow{r} \& \Tcal\\
\end{tikzcd}
\]

which shows that $\Lcal'$ does not depend on any choice of the presentation, and hence can be denoted $f^\#(\Lcal)$.
}

\subsection{Positivity and fiberwise density}

\block{\Def{\begin{itemize}
\item A locale $\Lcal$ is said to be positive, if whenever we can write $\Lcal$ as a union of open sublocales:

\[ \Lcal = \bigvee_{i \in I} u_i \]

the set of indices $I$ has to be inhabited. In this case, we write $\Lcal > \emptyset$.

\item A locale $\Lcal$ is said to be locally positive if every open sublocale can be written as a union of positive open sublocales.

\end{itemize}
}

If one assumes the law of excluded middle, then an open sublocale is positive if and and only if it is non-zero and every locale is locally positive (any non-zero element is the union of just itself, and the zero element is the empty union). But without the law of excluded middle this becomes a non trivial property.
 }
 
\block{If $X$ is a locale (preferably locally positive) we will denote by $\Ocal(X)^+$ the subset of positive open sublocales of $X$.}

\block{Local positivity is closely related to the notion of open map:

\Prop{Let $f: \Lcal \rightarrow \Mcal$ be a morphism of locale, then the following conditions are equivalent:

\begin{itemize}

\item For any $U$ open sublocale of $\Lcal$, its image $f_!(U)$ is an open sublocale of $\Mcal$; i.e. $f$ is an open map.

\item The frame morphism $f^* :\Ocal(\Mcal) \rightarrow \Ocal(\Lcal)$ has a left adjoint $f_{\circ}$ (i.e. $f_{\circ}(U) \leqslant V$ if and only if $U \leqslant f^*(V)$) which satisfies the additional identity:

\[ f_{\circ}( U \wedge f^*(V)) = \left( f_{\circ} U \right) \wedge V; \]

\item $\Lcal$ is locally positive as a $\Mcal$-locale.

\end{itemize}

Moreover in this situation, $f_{\circ}$ is the same as $f_!$ (restricted to open sublocales) and it corresponds to the internal map which associates to every $U \in \Ocal(\Lcal)$ the $\Mcal$-proposition `` $U$ is positive ".
}

For a proof, see \cite{borceux3}1.6.1 and 1.6.2 for the equivalence of the first two points, and see \cite{sketches} C3.1.17 for the last point.

}

\blockn{Because of this proposition, locally positive locales are generally called ``open locales". We cannot use this terminology here because we will have to speak a lot about locally positive sublocales, and ``open sublocales" would have two possible meaning in this case. The name ``overt" has also been proposed to avoid this confusion. }

\block{\label{locposlemma}The following lemma will often be useful to prove that some locales are locally positive:
\Lem{Let $X$ be a locale, and $p$ the morphism from $X$ to the point $\{*\}$. Assume that there is a basis $(b_i)_{i \in I}$ of $X$  and a collection of propositions $(w_i)_{i \in I}$ such that:

\[ w_i \Rightarrow (b_i)> \emptyset \]
\[ b_i \leqslant p^* w_i \]

Then $X$ is positive, $w_i$ is equivalent to $``b_i > \emptyset"$ and an arbitrary open sublocale of $X$ is positive if and only if it contains one of the $b_i$ such that $b_i > \emptyset$.
}

\Dem{As the $b_i$ form a basis, any $U \in \Ocal(X)$ can be written as:

\[ U = \bigvee_{i \in I \atop b_i \leqslant U} b_i \]

but as $b_i \leqslant p^*(w_i) = \bigvee_{w_i} \top $ one has:

\[ U = \bigvee_{i \in I \atop b_i \leqslant U} p^*(w_i) \wedge b_i = \bigvee_{i \in I \atop b_i \leqslant U \text{ and }w_i} b_i \]

as $w_i$ implies that $b_i$ is positive, this is an expression of $U$ as a supremum of positive open sublocales, proving that $X$ is locally positive. Now $w_i \Rightarrow b_i > \emptyset$ and as $b_i = \bigvee_{w_i} b_i$ one also has $b_i >\emptyset \Rightarrow w_i$, which proves the equivalence between $w_i$ and ``$b_i$ is positive". Finally if $U$ is positive, then from the previous expression of $U$ as a union, there exists an $i$ such that $b_i \leqslant U$ and $w_i$ hence $b_i$ is positive, and conversely if $U$ contains a positive $b_i$ then $U$ is itself positive.

}

}

\block{\label{prop_sitelocpos}\Prop{A locale $\Lcal$ is locally positive if and only if it can be defined by a Grothendieck site where each covering is inhabited. In this situation, an open $U$ of $\Lcal$ is positive if and only if it contains one of the representable.} This is essentially the localic version of \cite[C3.1.19]{sketches}. It can be applied to site as defined in \cite[C2.1.1]{sketches}, that is where the cover are only assumed to satisfies the base change axiom. }

\block{\label{pullbackopen}\Prop{Let $X$ be a locally positive locale in a topos $\Tcal$ and $f : \Ecal \rightarrow \Tcal$ a geometric morphism. Then $f^\#(X)$ is also locally positive, and (internally in $\Ecal$) an open $f^*(U) \in f^*(\Ocal(X))$ is positive if and only if $f^*(``U>\emptyset")$. }

\Dem{If one has a site of definition $(C,J)$ for $\Lcal$ in which each covering relation is inhabited then $f^*(C,J)$ also has this property and it is a site of definition for $f^\#(\Lcal)$. Hence this is an immediate corollary of the previous proposition.}
}

\block{\label{localicAC}Once we replace the idea of ``having points" by ``being positive and locally positive" to state that a locale is inhabited one can obtain a constructive version of ``the axiom of choice" in the form of:

\Prop{Let $I$ be a set with decidable equality and let $(X_i)_{i \in I}$ be a family of positive and locally positive locales. Then $\prod_{i \in I} X_i $ is positive and locally positive.
}

Note that the hypothesis that $I$ is decidable cannot be removed, and in fact cannot be weakened at least if we want to keep a first order property. See \cite[Chapter 3, 2.3.8]{mythesis} for more details. We are gratefull to Graham Manuell for pointing out a mistake in the original proof of this proposition.

\Dem{ For a finite product this follows from the fact that open surjections are stable under composition and pull-back (\cite[C3.1.11]{sketches}). If $I$ is decidable, then as a frame $\Ocal(\prod_{i \in I} X_i )$ is the directed colimits of the $\Ocal(\prod_{i \in P} X_i)$ for $P \subset I$ a finite subset, and all the transition map between the $\Ocal(\prod_{i \in P} X_i)$ are open surjections. We can then essentially copy the proof of \cite[C3.1.22]{sketches}: For each finite $P \subset I$ we have a site of definition for $\prod_{i \in P} X_i$ given by all the locally positive open sublocales and all covering relation between them. Then, because for $P \subset P' \subset I$ the transition map $\pi: \prod_{i \in P'} X_i \to \prod_{i \in P} X_i $ is an open surjection, the $\pi^*$ preserves locally positive elements and hence one can obtain a site of definition for $\Ocal(\prod_{i \in P} X_i)$ by taking the direct colimits of all these sites. Similarly to what happen in the proof of \cite[C3.1.22]{sketches}, the actual Grothendieck topology on this site is complicated to describe concretely: it is the smallest topology generated by union of the topologies on the sites for $\prod_{i \in P} X_i $. But the set of covers coming from all the $\prod_{i \in P} X_i $ already satisfies the base change axiom, that is, it is a site in the sense of Definition \cite[C2.1.1]{sketches}, so we can apply proposition \ref{prop_sitelocpos} to it and conclude that the product is indeed open.

  % This product can be explicitely describe by a site in which covering are all non-empty, which imply the result by proposition \ref{prop_sitelocpos}. Indeed we take as basic open of the product the open of the form

  % \[\bigcap_{i \in I} \pi_i^*(U_i) \]

  % Open surjections are stable by composition and pull-back (\cite[C3.1.11]{sketches}), hence if $X_1, \dots X_n$ are locally positive locales, then $\prod_{i=1}^n X_i$ also is. In general, a base of open sublocales of $\prod_{i \in I} X_i $ is given by the finite intersections of open sublocales of the form $\pi_i^*(U)$ for $U$ as open of $X_i$. If $I$ is decidable, each of these open sublocales can be rewritten as an intersection $\pi_{i_1}^*(U_1) \wedge \dots \wedge pi_{i_k}^*(U_k)$ with the $i_j$ pairwise distinct. Moreover, as each $X_i$ is locally positive such an open sublocale can be written as a union of sublocales of the same form but with the $U_i$ positive. As locales these open sublocales can be identified with $\prod_{i=1}^n U_i$ which is positive if each $U_i$ is positive, hence one has given a basis of positive open sublocales of the product $\prod_{i \in I} X_i $. This concludes the proof. 
}

}

\block{\label{localicADC}We also have a constructive version of the axiom of dependent choice:

\Prop{Let $X$ be an inhabited set equipped with a relation $R$ such that for each $x \in X$ there exists $ y \in X$ with $x R y$. Then the sublocale of $X^\N$ which classifies the sequences $(x_n)$ such that for each $n$ one has $x_n R x_{n+1}$ is positive and locally positive. }

This is proved in \cite{moerdijk1986connected} as lemma $C$.
}

\block{A geometric morphism $f :\Mcal \rightarrow \Lcal$ is said to be \emph{fiberwise dense} (or to have a fiberwise dense image) if for any proposition $U$, one has the relation:

\[ p^*(U) = f_*f^* p^* (U) \]

where $p$ denotes the canonical map $\Lcal \rightarrow \{*\}$ and $U$ is identified with an open sublocale of $\{*\}$.

A sublocale $S \subset \Lcal$ is said to be \emph{fiberwise closed} if it is fiberwise dense in no other sublocale of $\Lcal$.
}

\block{In the presence of the law of excluded middle these are equivalent to the more classical notions of density and closeness, but in general fiberwise density only implies density, and closeness only implies fiberwise closeness. For this reason they have also been called ``strongly dense" and ``weakly closed", but we prefer the terminology ``fiberwise" which is more uniform, more specific and allows less confusions. This name ``fiberwise" comes from the fact that, when interpreted internally in $\sh(X)$ for a (nice enough) topological space $X$, it indeed corresponds to a notion of fiberwise density (and fiberwise closeness) of morphisms of locales over $X$ whereas the usual notion of density would correspond to simple density, without taking the basis into account.

Aside from this difference of terminology, these definitions and the proof of all the results stated here can be found in \cite{sketches} after C1.1.22 and after C1.2.14.

Of course every sublocale $S$ admits a fiberwise closure $\overline{S}$ which is the smallest fiberwise closed sublocale containing $S$, or equivalently, the unique fiberwise closed sublocale in which $S$ is fiberwise dense.

}

\block{\label{fiberwisevsstrong}In the case of locally positive locales, the fiberwise density takes the following simpler form.
\Prop{Let $f:X \rightarrow Y$ be a map with $X$ locally positive. Then the following conditions are equivalent:

\begin{itemize}
\item[(a)] $f$ is fiberwise dense.

\item[(b)] $Y$ is locally positive, and for any positive open sublocale $U$ of $Y$, $f^*(U)$ is positive.
\end{itemize}

}

In presence of the law of excluded middle, every locale is locally positive and a positive open sublocale is just a non-zero open sublocale. Hence the previous proposition asserts (in presence of the law of excluded middle) that $f$ is fiberwise dense if for every non zero open sublocale $f^*(U)$ is also non zero, which is a classical characterisation of a dense map.

}

\block{\label{surjectimpopen}\Cor{Let $f:X \rightarrow Y$ be a surjection with $X$ locally positive, then $Y$ is locally positive. }
\Dem{A surjection is in particular a fiberwise dense map.}
}

\block{\label{fwdensesloc}\Prop{A fiberwise dense sublocale of a locally positive sublocale is also locally positive.}
}

\block{\label{pullbackfwdense}\Prop{If $g : X \rightarrow Y$ is a fiberwise dense map between two locally positive locales, then any pull-back of $g$ by a geometric morphism is also fiberwise dense. }

A counterexample to this proposition without the local positivity assumption can be found in \cite{sketches} right after corollary $C.1.2.16$.
}

\block{\Def{A locale $\Lcal$ is said to be weakly spatial if there exists a fiberwise dense map $P\rightarrow \Lcal$ with $P$ a spatial locale (or simply, with $P$ a set).}

By \ref{fiberwisevsstrong}, a weakly spatial locale is automatically locally positive, and a locally positive locale is weakly spatial if and only if every positive open sublocale has a point.

}

\block{\label{pullback-qdecid}\Lem{Let $X$ be any object of the base topos, then there exists a positive locally positive locale $\Lcal$, with $p$ the canonical geometric morphism from $\sh( \Lcal)$ to the base topos, such that $p^*X$ is the quotient of an object $I$ of $\sh(\Lcal)$ which has decidable equality. }

\Dem{One can take $\Lcal$ to be the classifying space for partial surjective maps from $\N$ to $X$. It is always a positive locally positive locale (see \cite{joyal1984extension}$V.3$ just after proposition $2$), and in $\sh(\Lcal)$ the object $p^* X$ is naturally a quotient of a subobject of $\N$, which is decidable. }

}

\block{\label{potentiallyweaklyspatial}\Prop{Let $X$ be a locally positive locale (of the base topos), then there exists a topos $\Tcal$ (even a locale) such that the canonical geometric morphism $p:\Tcal \rightarrow *$ is an open surjection and such that $p^\#(X)$ is weakly spatial in $\Tcal$. }

This result will be extremely important in the rest of this paper: indeed weak spatiality will play the same role as spatiality for complete metric spaces (see \ref{Metricset}), and as locales descend along open surjections this result will roughly allow us to assume whenever needed that all the metric locales involved come from metric sets.

\Dem{Thanks to the previous lemma, one can construct a locale $\Lcal$ in which one has a basis $(U_i)_{i \in I}$ of positive open sublocales of $p^\#(X)$ indexed by a set with decidable equality. By \ref{localicAC}:

\[ Y = \prod_{i \in I} U_i \]

is a positive locally positive locale, and corresponds to an open surjection (also denoted $p$) $p:\sh_{\Lcal}(Y) \rightarrow \Lcal \rightarrow * $. We will now prove that $p^\#(X)$ is weakly spatial. 

Internally in $\Lcal$, there is a canonical map $s_i : Y \rightarrow X \times Y$ defined as the composition of the i-th projection and the inclusion of $U_i$ into $X$ on the first component and the identity of $Y$ on the second component.  This defines a map of locale over $Y$:

\[ s :\coprod_{i \in I} Y \rightarrow X \times Y = p^\#(X) \]

which internally in $\sh_{\Lcal}(Y)$ gives a map $s$ from $f^*(I)$ to $p^\#(X)$ such that for each $i$, $s(i)$ is a point of $U_i$. As any positive open sublocale of $p^\#(X)$ contains one of the $U_i$, it shows that $p^\#(X)$ is weakly spatial.
}

}

\subsection{Descent theory} 

\label{Descenttheory}

\blockn{Let $\Ccal$ be a functor from the $2$-category of toposes to the $2$-category of categories,  like for example the functor which sends every topos $\Tcal$ to the category of internal locales of $\Tcal$, and any geometric morphism $f$ to the functor $f^{\sharp}$. We will denote by $f^*$ the action of a geometric morphism $f$ on $\Ccal$.}

\blockn{Let $f: \Ecal \rightarrow \Tcal$ be a geometric morphism, and let $c \in |\Ccal(\Ecal)|$. A descent data on $c$ is the data of an isomorphism $\epsilon : \pi_1^*(c) \rightarrow \pi_2^*(c) \in \Ccal(\Ecal \times_{\Tcal} \Ecal)$, such that if $\Delta$ denotes the diagonal map $\Delta: \Ecal \rightarrow \Ecal \times_{\Tcal} \Ecal$ then $\Delta^*( \epsilon)$ identifies with the identity map of $c$, and if $\pi_{1,2},\pi_{1,3}$ and $\pi_{2,3}$ denote the three projections $\Ecal \times_{\Tcal} \Ecal \times_{\Tcal} \Ecal \rightarrow \Ecal \times_{\Tcal} \Ecal$ and $\pi_1$, $\pi_2$ and $\pi_3$ the three projections from $\Ecal \times_{\Tcal} \Ecal \times_{\Tcal} \Ecal$ to $\Ecal$ then one has a commutative diagram:

\[
\begin{tikzcd}[ampersand replacement=\&]
\pi_1^*(c) \arrow{rd}{\pi_{13}^* \epsilon} \arrow{r}{\pi_{12}^*\epsilon} \& \pi_{2}^*(c) \arrow{d}{\pi_{23}^* \epsilon} \\
\& \pi_3^* (c)\\
\end{tikzcd}
\]

We define $Des(f,\Ccal)$ to be the category of objects of $\Ccal( \Ecal)$ endowed with a descent data (and morphisms being the morphisms in $\Ccal(\Ecal)$ whose pull-back along $\pi_1$ and $\pi_2$ commute to the $\epsilon$). If $c_0 \in \Ccal(\Tcal)$ then $f^* c$ is naturally endowed with a descent data and this defines a functor from $\Ccal(\Tcal)$ to $Des(f,\Ccal)$. One says that objects of $\Ccal$ descend along $f$, or that $f$ is a descent morphism\footnote{We follow the terminology of \cite{sketches}, it is in fact more common to say that $f$ is an effective descent morphism.} for $\Ccal$ if this functor induces an equivalence between $\Ccal(\Tcal)$ and $Des(f,\Ccal)$.

It is for example proved in \cite{joyal1984extension} that both objects and locales descend along open surjections. That is, for $\Ccal(\Tcal)= \Tcal$ and $\Ccal(\Tcal)$ being the category of internal locales of $\Tcal$ the geometric morphisms which are open and surjective are descent morphisms.
}

\blockn{In another language, the fact that objects of $\Ccal$ descend along all open surjections, or more generally along all geometric morphisms belonging to some Grothendieck topology one the $2$-category of topos exactly means that $\Ccal$ is a \emph{stack} for this topology. }

\subsection{Spaces of numbers}

\block{As mentioned in the introduction we are assuming that the base topos has a natural number object denoted by $\N$ (see \cite[A2.5 and D5.1]{sketches}). And from this natural number object one defines as usual the set $\Z$ of relative integers and $\Q$ of rational numbers with all their usual operations and properties.}

\block{$\R$ will denote the formal locale of real numbers, i.e;. classifying locale of the geometric propositional theory of Dedekind real numbers (continuous real number). When it is spatial (for example in presence of the law of excluded middle) it is the set of real numbers endowed with its classical topology. In any case, it agrees with the localic completion (as we define in \ref{compdef}) of $\Q$ for the Archimedean distance. $\C$ denote the formal locale of complex numbers, i.e. $\R \times \R$ endowed with its usual multiplication and addition.}

\block{Similarly will define a locale $\Length$ in which the distance function will take value. As earlier work of C.J.Mulvey showed we only care about knowing when a distance is smaller than some rational number, hence $\Length$ will be defined as the classifying locale of the theory of $P \subset \Q_+^*$ such that if $q \in P$ and $q<q'$ then $q' \in P$ and if $q \in P$ then there exists $q'<q$ such that $q' \in P$.

As $P$ is defined as a subset of positive rational numbers, $\Length$ corresponds only to non-negative numbers, and as we do not ask $P$ to be inhabited, $\Length$ contains a point $+\infty$ (corresponding to $P= \emptyset$). The topology on $\Length$ is the topology of upper semi-continuity i.e. the basic open sublocales are the $[0,q[$ for $q$ a rational (or real) number.
}

\block{On a topological space (or more generally in a Grothendieck topos) Dedekind real numbers correspond to continuous functions to $\R$, whereas points of $\Length$ correspond to non negative upper semi-continuous (possibly infinite) functions. This explains why Dedekind reals are called ``continuous" real numbers, and why points of $\Length$ can be called upper semi-continuous real numbers.}

\subsection{$[X,\R]$ is locally positive}
\blockn{The goal of this subsection is to show that, when $X$ is a compact regular locale, the locale $[X,\R]$ is locally positive (and hence also $[X,\C] \simeq [X,\R]^2$).}

\blockn{If $U$ and $V$ are two open sublocales of $X$ we write $U \ll V$ if $U$ is way below $V$, i.e. if when $V \leqslant \bigvee_{i \in I} U_i$ then there exists a finite subset $J \subset I$ such that $U \leqslant \bigvee_{j \in J} U_j$. We write $U \prec V$ when $U$ is rather below $V$, i.e. when $V \vee \neg U = X$, where $\neg U$ is the biggest open sublocale disjoint from $U$. A locale $X$ is regular when $\forall V \in \Ocal(X)$, $V = \bigvee_{U \prec V} U$. In a compact regular locale the two relations $\prec$ and $\ll$ are equivalent.
}

\blockn{In \cite{hyland1981function} one can find a description of the geometric theory classified by $[X,\R]$. This description shows that the open sublocales of the form $(U,q,q') = \{f | U \ll f^*(]q,q'[) \}$\footnote{Of course, we do not mean the set of points $f$ of $[X,\R]$ satisfying this properties, but the open sublocale classifying such functions $f$.} for $U$ an open sublocale of $X$ and $q$, $q'$ two rational numbers form a pre-basis of the topology of $[X,\R]$. 

\bigskip

As:

\[ U \ll f^*(]q,q'[)  \Leftrightarrow ( U \ll f^*(]q,+\infty[) ) \wedge( U \ll f^*(]-\infty,q'[) ), \]

$[X, \R]$ has a basis of open sublocales of the form

\begin{equation}\label{basicopen}
B = \left( \bigwedge_{i=1}^n (U_i,u_i,-) \right) \wedge \left( \bigwedge_{j=1}^m (V_j,v_j,+)  \right), 
\end{equation}

where $U_i$ and $V_i$ are open sublocales of $X$, $u_i$ and $v_i$ are rational numbers, $(U_i,u_i,-)$ denotes $\{f | U \ll f^*(]-\infty, u_i [) \}$ and $(V_j,v_j,+)$ denotes $\{f | V_j \ll f^*(]v_j,+\infty[) \}$.

}

\block{\Def{An open sublocale of the form given in $(\ref{basicopen})$ will be called a basic sublocale.
A basic sublocale will be said to be admissible if it satisfies the following condition:

\[ \forall i \in {1,\dots,n}, j \in {1,\dots, m}, ( u_i \leqslant v_j) \Rightarrow (\neg U_i) \vee (\neg V_j) = X.  \]
}

We will show in \ref{X,Clocpos} that a basic open is admissible if and only if it is positive, hence the property of being admissible is indeed a property of the open sublocale $B$, and not of its representation. But, while we have not proven this, we will assume that each time we consider a basic open $B$, it is given with a representation in the form of $(\ref{basicopen})$ and say that it is admissible if and only if its representation is.

}

\block{\label{lemmabasicpos} The following lemma is in some sense a constructive form of Urysohn's lemma, asserting that compact regular locales are in fact completely regular.
\Lem{Let $X$ be a compact regular locale, and let $U$,$V$ be two open sublocales of $X$ such that $U \ll V$. Then there exists a positive locally positive locale $\Lcal$, such that in the logic of $\Lcal$ there exists a continuous function from $X$ to $[0,1]$\footnote{That is externally a function from $\Lcal \times X$ to $[0,1]$.} such that $f$ restricted to $U$ is zero and $f$ is constant equal to one on $\neg V$. }

\Dem{The classical proof of the Urysohn lemma for locale (see for example \cite[Chap. XIV]{picado2012frames}) goes as follows: In a compact regular locale the relation $U \prec V $ is equivalent to the relation $U \ll V$. The relation $\prec $ in general does not interpolate, but in a locally compact locale the relation $ \ll $ always does, ie if $a \ll b$ then there exists $c$ such that $a \ll c \ll b$. In particular in a compact regular space the relation $\prec$ interpolates and (using the axiom of choice) one can construct a $\Q$-indexed family of open subspaces $U_q$ such that $U_0 = U$ , $U_1 = V$ and if $q <q'$ then $U_q \prec V_{q'}$, and we define $U_q = \emptyset$ when $q<0$ and $U_q = X$ when $q>1$. This defines a ``scale" (see \cite{picado2012frames} XIV.5.2 ) which in turns defines a function from $X$ to $[0,1]$ with the required property (see \cite{picado2012frames}XIV.5.2.2).

\bigskip

The only part of the previous proof which is not constructive is the application of the axiom of dependent choice to construct the sequence $U_q$. By applying \ref{localicADC} one can construct a locale $\Lcal$ in which there exists such a sequence and then finish the proof in the logic of $\Lcal$ by constructing the function we are looking for. The only thing we need to check is that if $x \prec y$ then their pull-back to $\Lcal$ also satisfy this identity, but as it can equivalently be defined by  `` $\exists c$ such that $x \wedge c = \emptyset$ and $c \vee y = \top$ " this is immediate.

}

}

\block{\label{admihavepoints}
\Prop{If $X$ is compact completely regular and $B$ is an admissible basic sublocale of $[X,\R]$, then $B$ has a point. If $X$ is just compact regular and $B$ is admissible then $B$ is positive.}
\Dem{Assume that $X$ is completely regular, and let us first remark that when $X$ is a compact completely regular locale, if $U$ and $V$ are two open sublocales of $X$ such that $(\neg U) \vee (\neg V) =X$, then, as $U \ll (\neg V)$, it is possible to construct a continuous function $f:X \rightarrow [0,1]$ such that $f$ restricted to $U$ is constant equal to $0$ and $f$ restricted to $V \subseteq \neg \neg V $ is constant equal to $1$.

\bigskip

Now let
\[ B = \left( \bigwedge_{i=1}^n (U_i,u_i,-) \right) \wedge \left( \bigwedge_{j=1}^m (V_j,v_j,+)  \right) \]
be an admissible basic sublocale of $[X,\R]$.

Let $\epsilon$ be a positive rational number smaller than all the positive differences between two numbers of the form $u_i$ or $v_i$. For each couple $(i,j)$ we choose a continuous function $f_{i,j} : X \rightarrow \R$ such that:

\begin{itemize}

\item If $v_j < u_i$ then $f_{i,j}$ is the constant function equal to $\frac{v_j+u_i}{2}$

\item If $u_i \leqslant v_j$  then  $(\neg U_i) \vee (\neg V_j) =X$ and $f_{i,j}$ is a continuous function such that $f$ is constant equal to $u_i-\epsilon$ on $U_i$, $f$ is constant equal to $v_j+\epsilon$ on $V_j$ and $f$ takes value in $[u_i-\epsilon,v_j+\epsilon]$. (such a function exists by the previous remark).

\end{itemize}

Then,
\[ f = \max_{1 \leqslant j \leqslant m } \min_{1 \leqslant i \leqslant n} f_{i,j}, \]

is a point of $B$. Indeed:

\begin{itemize}
\item Let $i \in \{1,\dots,n\}$, then (on $U_i$), since for each $j$, $f_{i,j}$ is smaller than $u_i-\frac{\epsilon}{2}$, the infimum $\inf_{i'=1}^n f_{i',j}$ is smaller than $u_i-\frac{\epsilon}{2}$ and $f$ smaller than $u_i-\frac{\epsilon}{2}$ on $U_i$ as a (finite) supremum of a quantities smaller than  $u_i-\frac{\epsilon}{2}$.

\item Let $j \in \{1,\dots,m\}$, then (on $V_j$), as for each $i$, $f_{i,j}$ is greater than $v_j + \frac{\epsilon}{2}$, the infimum $\inf_{i=1}^n f_{i,j}$ is greater than $v_j + \frac{\epsilon}{2}$. And $f$ is greater than $v_j + \frac{\epsilon}{2}$ on $V_j$.

\end{itemize}

This concludes the proof when $X$ is completely regular. We now assume that $X$ is only regular. Then all the functions $f_{i,j}$ we used in the first part can be instead constructed in the logic of positive locally positive locales $\Lcal_{i,j}$ using \ref{lemmabasicpos}. The product $\Lcal$ of all these $\Lcal_{i,j}$ is also positive and locally positive by \ref{localicAC}, and in the logic of $\Lcal$, all the functions $f_{i,j}$ we used in the first part exist and hence one can construct the function $f$ which is going to be a point of $B$ in the logic of $\Lcal$ exactly as we did above. This defines a map $\Lcal \rightarrow B $ and, as $\Lcal$ is positive, this proves that $B$ is positive and concludes the proof. }

}

\block{\label{posimpadmis}\Lem{Let $p$ denote the canonical map from $[X,\R]$ to the point. Let $B$ be a basic sublocale then one has:

\[ B \leqslant p^*( \text{`` $B$ is admissible " })  \]

where we identify the proposition ``$B$ is admissible" with a subset of $\{*\}$ and hence with an open sublocale of the point.
}

\Dem{ We will prove that in the theory classified by $[X, \R ]$ (describe in \cite{hyland1981function}) the proposition asserting that $B$ is admissible can be deduced from the proposition corresponding to $B$. 

Indeed, let $B$ be as in $(\ref{basicopen})$ and let $i$ and $j$ such that $u_i \leqslant v_j$.

one has:

\[ B \vdash \left( U_i \ll f^*(]-\infty, u_i[)\right) \wedge \left( V_j \ll f^*(]v_j,+\infty[) \right), \]

\[ \left( U_i \ll f^*(]-\infty, u_i[)\right) \vdash \bigvee_{U_i \ll U } \left( U \ll f^*(]-\infty, u_i[)\right) \]

and 

\[ \left( U \ll f^*(]-\infty, u_i[)\right) \wedge \left( V \ll f^*(]v_j,+\infty[) \right) \vdash (U \wedge V)= \emptyset.  \]

Hence

\[ B \vdash \bigvee_{U_i \ll U \atop V_j \ll V } (V \wedge U = \emptyset) \]

but for any $U_i \ll U$ and $V_j \ll V$ if $(V \wedge U = \emptyset)$ then $\neg U \vee \neg V =X$ because

\[ 
\begin{array}{c c l}
 X& = &(\neg U_i \vee U) \wedge (\neg V_j \vee V)  \\ &=& (\neg U_i \wedge \neg V_j) \vee ( \neg U_i \wedge V) \vee (U \wedge \neg V_j) \vee ( U \wedge V) \end{array}
\]
The last term of the union can be removed by assumption, and we can duplicate the first, obtaining

\[ \begin{array}{c c l}
 X & = & [( \neg U_i \wedge \neg V_j) \vee ( \neg U_i \wedge V)] \vee [(U \wedge \neg V_j) \vee ( \neg U_i \wedge \neg V_j)] \\
 &=& [(\neg U_i) \wedge (\neg V_j \vee  V)] \vee [(\neg V_j) \wedge (\neg U_i \vee U)] \\
 &=& \neg U_i \vee \neg V_j 
 \end{array}\]
 
Hence $B \vdash   \neg U_i \vee \neg V_j$. As this is true for any $(i,j)$ such that $u_i \leqslant v_j$ we get the desired result.

}

}

\block{\label{X,Clocpos} Combining all these results we obtain:

\Th{If $X$ is a compact regular locale, then a basic sublocale $B$ of $[X,\R]$, is admissible if and only it is positive. In particular, $[X,\R]$ is locally positive and the admissible basic sublocales form a basis of positive open sublocales. }

\Dem{It suffices to apply Lemma \ref{locposlemma} with $b_i$ the basic open sublocales and $w_i$ the propositions ``$b_i$ is admissible". Proposition \ref{admihavepoints} shows that $w_i$ implies $b_i > \emptyset$ and \ref{posimpadmis} is exactly the second condition.
}
} 

\block{\label{compregandwspa} We also obtain the following

\Prop{Let $X$ be a compact regular locale, $X$ is completely regular if and only if $[X, \R ]$ is weakly spatial.}

\Dem{If $X$ is completely regular, then \ref{admihavepoints} shows that each admissible has a point. But by \ref{X,Clocpos} they form a basis of positive open, hence this proves that points of $[X, \R]$ are dense. Conversely, if $[X,\R]$ is weakly spatial and $U$,$V$ are two open sublocales of $X$ such that $U \prec V$, then there exists $W$ such that $U \prec W \prec V$ and the basic open:

\[ B = (U,0,-) \wedge (\neg W,1,+) \]

is admissible because $\neg U \vee \neg \neg W \geqslant \neg U \vee W =X$. Hence it is positive and hence it has a point. But a point of $B$ is a function from $X$ to $\R$ such that $f$ is negative on $U$ and greater than one on $\neg W$. As $\neg W \vee V =X$ the function $f$ shows that $U$ is ``completely below $V$", and this proves that $X$ is completely regular.
}

 }

\section{Constructive theory of metric locales}
\label{sectionmetriclocale}
\subsection{Pre-metric locale}
\blockn{As our major concern is the study of localic Banach spaces, we will only consider metrics on a locale which are defined by a distance function. However, it should be noted that the point \ref{PMP_distfromdiam} of the series of propositions given in \ref{premetricprop} shows that one can specify a distance by giving the diameter $\delta(U)$ of each open sublocale $U$, and the classical theory\footnote{Which has not been done constructively yet as far the author knows.} which can be found for example in the chapter XI of \cite{picado2012frames} suggests that a definition by diameters should also be possible.}

\block{\Def{A pre-distance $d$ on a locale $X$ is a function 

\[ d: X \times X \rightarrow \Length \]

which is symmetric ($d(x,y)=d(y,x)$), satisfies the triangular inequality $d(x,y) \leqslant d(x,z) + d(z,y)$ and such that $d(x,x)=0$

A pre-metric locale is a locally positive locale $X$ endowed with a pre-distance.
}

We insist on the fact that our pre-metric locale are always assumed to be locally positive. We do not know exactly which parts of the theory of metric locales it is possible to develop without this hypothesis (without it, one should at least avoid everything which uses the construction $B_q \Lcal$ of \ref{premetricdef} but it seems that what is left is relatively well behaved without it). In any case, the theory is at least easier, and probably nicer with this local positivity assumption. Theorem \ref{X,Clocpos} shows that this case is enough for the Gelfand duality, and as locale positivity descend along open surjections and is automatic for metric sets it is also enough to obtain good descent properties.

\smallskip

Of course, the formulas $d(x,y)=d(y,x)$ and $d(x,y) \leqslant d(x,z) + d(z,y)$ have to be interpreted in a diagrammatic way or in terms of generalized points. In particular, if we define
 \[\Delta_q := \{ (x,y) | d(x,y) < q \} = d^*\left(\overleftarrow{[0,q[}\right)\] 
then the symmetry means that $\Delta_q$ is invariant by exchange of the two factors, $d(x,x)=0$ means that for all $q$, $\Delta_q$ contains the diagonal embeddings of $X$, and finally the triangular inequality means that: 

\[ \pi_{1,2}^*(\Delta_q) \wedge \pi_{2,3}^*(\Delta_{q'}) \leqslant \pi_{1,3}^*(\Delta_{q+q'}) \]

Where $\pi_{i,j}$ denote the various projections from $X^3$ to $X^2$.
}

\block{\label{premetricdef}\Def{Let $X$ be a pre-metric locale, and $\Lcal$ and $\Mcal$ be two sublocales of $X$. then
\begin{itemize}
\item We say that $\delta(\Lcal)<q$ if $\Lcal \times \Lcal \subseteq \Delta_{q'}$ for some positive rational number $q'<q$. One easily sees that $\delta(\Lcal)$ is indeed an element of $\Length$;

\item We say that $\Lcal \triangleleft_q \Mcal$ if $\pi_1^*(\Lcal) \wedge \Delta_q \leqslant \pi_2^*(\Mcal)$. We say that $\Lcal \triangleleft \Mcal$ if $\Lcal \triangleleft_q \Mcal$ for some positive rational $q$;

\item if $q$ is a positive rational number then $B_q \Lcal = (\pi_2)_! (\pi_1^*(\Lcal) \wedge \Delta_q )$.

\end{itemize}
}
These should be interpreted as: $\delta$ is the diameter of a sublocale, $B_q$ is the $q$ neighborhood of a sublocale and $\Lcal \triangleleft_q \Mcal$ means that the $q$ neighborhood of $\Lcal$ is included in $\Mcal$.}

\block{We will denote by $\Ocal(X)^{<q}$ the set of open sublocales $U$ of $X$ such that $\delta(U)<q$, and $\Ocal(X)^{+,<q}$ will be simply the subset $\Ocal(X)^+ \cap \Ocal(X)^{<q}$ of positive elements of $\Ocal(X)^{<q}$.} 

\block{\label{premetricprop}\Prop{
\begin{enumerate}
\item $B_q \Lcal \subseteq \Mcal $ if and only if $ \Lcal \triangleleft_q \Mcal$.
\item If $\Lcal \subseteq \Mcal$ then $\delta(\Lcal) \leqslant \delta(\Mcal)$.
\item If $\Lcal \triangleleft \Mcal$ then $\Lcal \subseteq \Mcal$. In particular for all positive rational numbers $q$ one has $\Lcal \subseteq B_q \Lcal$.
\item \label{PMP_compatriang}If $\Lcal \triangleleft_q \Mcal$ and $\Lcal' \triangleleft_q \Mcal'$ then $\Lcal \wedge \Lcal' \triangleleft_q \Mcal \wedge \Mcal'$ and $\Lcal \vee \Lcal' \triangleleft_q \Mcal \vee \Mcal'$.

\item \label{PMP_diamassup}$\displaystyle \delta\left( \bigvee_{i \in I} \Lcal_i \right) = \sup_{i,j \in I} \delta(\Lcal_i \vee \Lcal_j)$

\item \label{PMP_triangulardiamaters} If $\Lcal \wedge \Mcal$ contains a positive and locally positive sublocale then $\delta(\Lcal \vee \Mcal) \leqslant \delta(\Lcal)+\delta(\Mcal)$.

\item \label{PMP_triangulardiamatersN} Let $(\Lcal_i)_{i=0 \dots n}$ be a finite sequence of sublocales such that for all $i$, $\Lcal_{i-1} \wedge \Lcal_{i}$ contains a positive and locally positive sublocale then:

\[\delta\left( \bigvee_{i=0}^n \Lcal_i \right) \leqslant \sum_{i=0}^n \delta(\Lcal_i) \]

\item For any $q>0$, $\Ocal(X)^{<q}$ is a basis of the topology of $X$.

\item \label{PMP_distfromdiam}$\displaystyle \Delta_q = \bigvee_{U \in \Ocal(X)^{<q}} U \times U $

\item \label{PMP_BqLunion} If $\Lcal$ is locally positive, then
\[\displaystyle B_q \Lcal = \bigvee_{U \in \Ocal(X)^{<q} \atop U \wedge \Lcal >\emptyset } U.\]

In particular, if $\Lcal$ is locally positive, $B_q \Lcal$ is open.

\item \label{PMP_triangularB}If $\Lcal$ is locally positive then
\[ B_{q'}(B_{q}(\Lcal)) \subseteq B_{q+q'} (\Lcal). \]

\item \label{PMP_diameterB} If $\Lcal$ is locally positive then $\delta(B_q \Lcal) \leqslant 2q+\delta(\Lcal)$.

\end{enumerate}
}
\Dem{
\begin{enumerate}
\item This is simply the adjunction between $(\pi_2)_!$ and $(\pi_2)^*$.
\item If $\Lcal \subseteq \Mcal$ and if $\delta(\Mcal) <q$ then there exists a positive rational $q'<q$ such that $\Lcal \times \Lcal \subseteq \Mcal \times \Mcal \subseteq \Delta_{q'}$ hence $\delta(\Lcal)<q$.
\item Assume that $\pi_1^*( \Lcal) \wedge \Delta_q \subseteq \pi_2 ^*(\Mcal)$ for some positive rational number $q$, and let $i: X \rightarrow X \times X$ be the diagonal embedding, then:

\[ i^*(\pi_1^*( \Lcal) \wedge \Delta_q ) \subseteq i^* \pi_2 ^*(\Mcal) = \Mcal \]

And:

\[ i^*(\pi_1^*( \Lcal) \wedge \Delta_q ) = i^* \pi_1^*(\Lcal) \wedge i^* \Delta_q =\Lcal \wedge X = \Lcal \]

hence $\Lcal \subseteq \Mcal$. The second part of the result then follows from the fact that as $B_q \Lcal \subseteq B_q \Lcal$, one has $\Lcal \triangleleft_q B_q \Lcal$.

\item Assume that $\pi_1^* \Lcal \wedge \Delta_q \subseteq \pi_2^* \Mcal$ and that $\pi_1^* \Lcal' \wedge \Delta_q \subseteq \pi_2^* \Mcal'$, then:

\[ \pi_1^*( \Lcal \wedge \Lcal') \wedge \Delta_q = \pi_1^*(\Lcal) \wedge \Delta_q \wedge \pi_1^*(\Lcal') \wedge \Delta_q \subseteq \pi_2^*(\Mcal) \wedge \pi_2^*(\Mcal') \]

hence $\Lcal \wedge \Lcal \triangleleft_q \Mcal \wedge \Mcal $.

And for the union:

\[ \begin{array}{r c l}
 \pi_1^*(\Lcal \vee \Lcal') \wedge \Delta_q &=& (\pi_1^*(\Lcal) \vee \pi_1^*(\Lcal') ) \wedge \Delta_q  \\ &=&( \pi_1^* \Lcal \wedge \Delta_q) \vee (\pi_1^* \Lcal' \wedge \Delta_q)\\ &\subseteq& \pi_2^*(\Mcal)\vee \pi_2^*(\Mcal'),
\end{array} \]

which gives the result.

The fact that intersections distribute over finite unions of sublocales and that pull-backs preserve finite unions of sublocales can be found in \cite{sketches} C1.1.15 and C.1.19, but formulated in terms of frames instead of locales (i.e. union of sublocales correspond to intersection of nuclei, and pull-back of a sublocale to a pushout).

\item Clearly, $\sup_{i,j \in I} \delta(\Lcal_i \vee \Lcal_j) \leqslant \delta\left( \bigvee_i \Lcal_i \right)$ because $\Lcal_i \vee \Lcal_j \subseteq \bigvee \Lcal_i$. Let $q$ such that $\sup_{i,j \in I} \delta(\Lcal_i \vee \Lcal_j)<q$ i.e. there exists $q'<q$ such that for all $i,j$, $\delta(\Lcal_i \vee \Lcal_j)<q'$. But as

\[ \left( \bigvee_{i \in I} \Lcal_i \right) \times \left( \bigvee_{j \in I} \Lcal_j \right) = \bigvee_{i,j} \Lcal_i \times \Lcal_j\]

and for all $i,j$, $\Lcal_i \times \Lcal_j \subseteq \Delta_{q'}$, one obtains

\[ \left( \bigvee_{i \in I} \Lcal_i \right) \times \left( \bigvee_{j \in J} \Lcal_j \right) \subseteq \Delta_{q'}, \]

which concludes the proof.

\item Assume that $\Lcal \times \Lcal \subseteq \Delta_{q}$ and $\Mcal \times \Mcal \subseteq \Delta_{q'}$, we will prove that, under the assumption of the proposition, $(\Lcal \vee \Mcal) \times (\Lcal \vee \Mcal) \subseteq \Delta_{q+q'}$.

As $(\Lcal \vee \Mcal) \times (\Lcal \vee \Mcal) = (\Lcal \times \Lcal) \vee (\Lcal \times \Mcal) \vee (\Lcal \times \Mcal) \vee (\Mcal \times \Mcal)$ and $(\Lcal \times \Lcal)$ and $(\Mcal \times \Mcal)$ are already known to be subsets of $\Delta_{q+q'}$, we only have to prove it for $(\Lcal \times \Mcal)$ and $(\Mcal \times \Lcal)$. In $X^3$ one has:

\[ \begin{array}{r c l} \Mcal \times (\Lcal \wedge \Mcal) \times \Lcal \subseteq \pi_{1,2}^*(\Mcal \times \Mcal) \wedge \pi_{2,3}^*(\Lcal \times \Lcal) &\subseteq& \pi_{1,2}^*(\Delta_q') \wedge \pi_{2,3}^*(\Delta_q) \\
& \subseteq & \pi_{1,3}^*(\Delta_{q'+q})
\end{array} \]

Applying $(\pi_{1,3})_!$ yields the result because as $(\Lcal \times \Mcal)$ contains some positive and locally positive sublocale, the projection $\pi_{1,3}$ from $\Lcal \times (\Lcal \wedge \Mcal) \times \Mcal$ to $\Lcal \times \Mcal$ is a surjection.

\item It is immediate by induction on $n$ using the previous point.

\item Thanks to the point $2.$ it is enough to check that $\Ocal(X)^{<q}$ covers $X$. Take a covering of $\Delta_{q/2}$ by open sublocales of the form $U_i \times V_i$, then pulling back along the diagonal embeddings of $X$ into $\Delta_{q/2}$ one has:

\[ X =\bigvee_i U_i \wedge V_i \]

but $(U_i \wedge V_i) ^2 \leqslant U_i \times V_i \leqslant \Delta_{q/2}$ hence $\delta(U_i \wedge V_i )<q$ which concludes the proof.

\item Thanks to the previous point, for any $q'<q$, $\Delta_{q'}$ can be written as a union of $U_i \times V_i$ with $\delta(U_i)<{q'}$ and $\delta(V_i)<{q'}$.
If $U_i \times V_i \subseteq \Delta_{q'}$. then so does $V_i \times U_i$, and hence, in our situation:

\[ (U_i \cup V_i)^2 = (U_i \times U_i) \cup (V_i \times U_i) \cup (U_i \times V_i) \cup (V_i \times V_i) \subseteq \Delta_{q'} \]

Hence $\delta(U_i \cup V_i)<q$ and the $(U_i \cup V_i)^2$ cover $\Delta_{q'}$. This being done for an arbitrary $q'<q$, these open sublocales also cover $\Delta_q$, because as the $\Delta_q$ are defined by a function from $X \times X$ to $\Length$ one has

\[ \Delta_q = \bigvee_{q'<q} \Delta_{q'} \]

\item Applying the definition of $B_q V$ using that $\pi_1^*(\Lcal) = \Lcal \times X$ and the previous point gives directly

\[B_q \Lcal = (\pi_2)_! \left( \bigvee_{\delta(U) < q }(\Lcal \wedge U) \times U \right) = \bigvee_{\delta(U)<q \atop \Lcal \wedge U > \emptyset} U. \]

\item From the previous point

\[ B_q(B_{q'} \Lcal) = \bigvee_{v \in \Ocal(X)^{<q} \atop v \wedge B_{q'} \Lcal > \emptyset} v \]

But, still by the previous point, an open sublocale $v$ of $X$ satisfies $v \wedge B_{q'} \Lcal > \emptyset$ if and only if there exists $v' \in \Ocal(X)^{<q'}$ such that $v' \wedge \Lcal >\emptyset$ and $v \wedge v'>\emptyset$. For any open sublocale of this sort, one has $\delta(v \vee v') < q+q'$ by point \ref{PMP_triangulardiamaters}. Hence $v \vee v'$ is a positive open sublocale such that $\delta(v \vee v') <q+q'$ and $(v \vee v') \wedge \Lcal > \emptyset$. In particular $v \leqslant v \vee v' \leqslant B_{q+q'} \Lcal $.

This proves that $B_q(B_{q'} \Lcal) \leqslant B_{q+q'} \Lcal$.

\item  From point \ref{PMP_BqLunion} one has
\[B_q \Lcal = \bigvee_{v \in \Ocal(X)^{<q} \atop v \wedge \Lcal >\emptyset} v. \]
Hence from point \ref{PMP_diamassup} one has
\[ \delta(B_q \Lcal) = \sup_{v,v' \in \Ocal(X)^{<q} \atop v \wedge \Lcal, v' \wedge \Lcal >\emptyset} \delta(v \vee v').\]

But for any two such $v,v'$ one has by point \ref{PMP_triangulardiamatersN}: $\delta(v \vee v' )\leqslant \delta(v \vee v' \vee \Lcal) \leqslant \delta(\Lcal) + \delta(v) + \delta(v') \leqslant \delta(\Lcal)+2q$. One obtains the result by taking the supremum.

\end{enumerate}

}

}

\block{\label{overlinedelta}Usually, the distance function $d:X \times X \rightarrow \Length$ is expected to be in fact a continuous map from $X \times X$ to $\mathbb{R}$, and not only a semi-continuous map as our definition of distance suggest it. The reason for our choice is that we know (see for example \cite{burden1979banach}) that the norm on a Banach space has to take value in $\Length$, even if we want to think of it as a function which is continuous\footnote{as opposed to semi-continuous.}. Classically, the continuity is a consequence of the triangular inequality, and the following proposition gives a constructive interpretation of this result, restoring a form of ``fiberwise continuity" of $d$.

\Prop{Let $\overline{\Delta_q}$ be the fiberwise closure of $\Delta_q$ in $X \times X$. Then for all $q<q'$ one has $\overline{\Delta_q} \subseteq \Delta_{q'}$.}

\Dem{Let $q'$ be a rational such that $q<q'$ and let $\epsilon = \frac{q'-q}{2}$. As $\Delta_q$ is by definition fiberwise dense in $\overline{\Delta_q}$, Proposition \ref{fiberwisevsstrong} implies that $\overline{\Delta_q}$ is locally positive, and in particular one can write that

\[ \overline{\Delta_q} \leqslant \bigvee_{v,v' \in \Ocal(X)^{<\epsilon} \atop v \times v' \wedge \overline{\Delta_q}> \emptyset } v \times v'. \]

But, still by \ref{fiberwisevsstrong} and by fiberwise density of $\Delta_q$ in $\overline{\Delta_q}$, for any two such $v,v'$ one has $v \times v' \wedge \Delta_q >\emptyset$ and hence there exists $U$ such that $\delta(U)<q$ and $(v \times v') \wedge (U \times U)$ is positive. This implies that $v \wedge U$ and $v' \wedge U$ are positive and hence, by point \ref{PMP_triangulardiamatersN} of \ref{premetricprop}, that $\delta(v \vee v') \leqslant \delta(v)+\delta(v')+\Delta(U)<q+2\epsilon =q'$.

Therefore,
\[ v \times v' \subseteq (v \vee v' ) \times (v \vee v') \subseteq \Delta_{q'}, \]

and this concludes the proof.

}
}

\block{\Def{Let $X$ be a pre-metric locale, we will say that $X$ has a continuous distance if the pre-distance function $d : X \times X \rightarrow \Length$ internally corresponds to a continuous real number, i.e. if the pre-distance function factors into $X \times X \rightarrow \overline{\R+} \rightarrow \Length$. In this situation we define $\Theta_q$ to be the open sublocale of $X \times X$ corresponding to $\{(x,y) | d(x,y)>q \}$.} }

\block{\label{LEMimpdistcont}Assuming the law of excluded middle, we indeed obtain continuity:
\Prop{Assuming the law of excluded middle in the base topos, any pre-metric locale has a continuous distance.} 

\Dem{If one assumes the law of excluded middle in the base topos then any fiberwise closed sublocale is in fact a closed sublocale. In particular, there exists open sublocales $\Theta'_q$ of $X \times X$, which are the complementary open sublocales of the (closed) sublocales $\overline{\Delta_q}$. From the fact, proved in \ref{overlinedelta} that for any $q<q'$ one has the relation

\[ \Delta_q \leqslant \overline{\Delta_q} \leqslant \Delta_{q'} \]

and we deduce

\[ \Delta_q \wedge \Theta'_q = \emptyset \]
\[ \Delta_{q'} \vee \Theta'_q = X \times X \]

and $\overline{\Delta_q} \leqslant \overline{\Delta_{q'}}$ gives $\Theta'_{q} \geqslant \Theta'_{q'}$.

If we define, $\Theta_q= \bigvee_{q<q'} \Theta'_{q'}$, then $\Delta_q$ and $\Theta_q$ define a map from $X \times X$ to $\overline{\R+}$ which yields the desired factorisation.
}
}

\block{\label{metricmap}\Prop{Let $f:X \rightarrow Y$ be a map between two pre-metric locales. Then the following conditions are equivalent:
\begin{enumerate}
\item[(a)] For any positive rational $q$, $ \Delta_q \subseteq (f \times f)^*(\Delta_q) $

\item[(b)] For any locally positive sublocale $\Lcal$ of $X$, $\delta(f_!\Lcal) \leqslant \delta(\Lcal)$.

\item[(c)] For any $U\in \Ocal(X)^{<q_1}$, $v_1 \in \Ocal(Y)^{<q_2}$, $v_2 \in \Ocal(Y)^{<q_3}$ such that $f^*(v_1) \wedge U$ and $f^*(v_2) \wedge U$ are positive, one has $\delta(v_1 \vee v_2) < q_1+q_2+q_3$.

\item[(d)] For any $U \in \Ocal(X)$ and any positive rational $q$: \[ \delta(B_q f_! U) \leqslant \delta(U) + 2 q. \]

\item[(e)] For any open sublocale $U$ of $X$ such that $\delta(U)<q$ there exists an open sublocale $V$ of $Y$ such that $\delta(V)<q$ and $U \subseteq f^*(V)$. 
\end{enumerate}

A map satisfying these conditions is called a metric map.
}

Of course, condition $(a)$  is the point free formulation of the usual $d(f(x),f(y))\leqslant d(x,y)$.

\Dem{
\begin{enumerate}

\item[$(a) \Rightarrow (b)$] Let $q$ such that $\delta(\Lcal)<q$, i.e. there exists $q'<q$ such that $\Lcal \times \Lcal \subseteq \Delta_{q'}$. Hence,

\[ \Lcal \times \Lcal \subseteq (f \times f)^* (\Delta_{q'}) \]

This proves that the image $(f \times f)_!(\Lcal \times \Lcal)$ in $X \times X$ is included in $\Delta_{q'}$. Unfortunately, as a product of surjections may fail to be a surjection, it is not enough to conclude directly that $f_!(\Lcal) \times f_!(\Lcal) \subseteq \Delta_{q'}$. But we can still conclude using the fact that as $\Lcal$ and $f_!(\Lcal)$ are both locally positive, then by \ref{pullbackfwdense} (applied twice) the map $f : \Lcal \times \Lcal \rightarrow f_!(\Lcal) \times f_!(\Lcal)$ is always fiberwise dense. This implies that $\Delta_{q'}$ is fiberwise dense in $f_!(\Lcal) \times f_!(\Lcal)$, and by \ref{overlinedelta} that:

\[ f_!(\Lcal) \times f_!(\Lcal) \subseteq \overline{\Delta_{q'}} \subseteq \Delta_{q}\]

which concludes the proof.

\item[$(b) \Rightarrow (c)$] by \ref{surjectimpopen} $\Lcal = f_!(U)$ is locally positive  because $U$ is and $f:U \rightarrow f_!(U)$ is a surjection. Also, $\delta(f_!(U))<q_1$ by $(b)$. Hence one obtains $(c)$ by applying point \ref{PMP_triangulardiamatersN} of \ref{premetricprop} (with n=2), together with the fact that $f^*v \wedge U >\emptyset$ is equivalent to $v \wedge f_!U >\emptyset$ because $f: U \rightarrow f_! U$ is a surjection and hence in particular a fiberwise dense map.

\item[$(c) \Rightarrow (d)$]One has

\[B_q f_! U = \bigvee_{v \in \Ocal(Y)^{<q} \atop f^*(v) \wedge U>\emptyset} v \]

The same argument as given for point \ref{PMP_diameterB} of \ref{premetricprop} allow one to conclude.

\item[$(d) \Rightarrow (e)$] If $\delta(U)<q$ then there exists a positive $\epsilon$ such that $\delta(U)<q - 2 \epsilon$. Take $V = B_\epsilon f_! U$ yields the result as $U \leqslant f^*f_! U \leqslant f^* B_{\epsilon} f_! U = f^* V$.

\item[$(e) \Rightarrow (a)$] 
Using $(e)$ one gets immediately the inclusion
\[ \Delta_q = \bigvee_{U \in \Ocal(X)^{<q}} U \times U \subseteq \bigvee_{V \in \Ocal(Y)^{<q}} f^*(V) \times f^*(V) = (f \times f)^*(\Delta_q) \]

\end{enumerate}
}
}

\block{\label{uniformmap}\Prop{Let $f:X \rightarrow Y$ be a map between two pre-metric locales, let $\epsilon$ and $\eta$ be two positive rational numbers, then the following conditions are equivalent:

\begin{enumerate}
\item[(a)] $\Delta_{\eta} \leqslant (f \times f)^* \Delta_{\epsilon}$
\item[(b)] If $U \in \Ocal(X)$ and $\delta(U)<\eta$ then $\delta(f_!(U))<\epsilon$
\item[(c)] If $U \in \Ocal(X)$ and $\delta(U)<\eta$ then there exists $V \in \Ocal(Y)$ such that $\delta(V)<\epsilon$ and $U \leqslant f^*(V)$.

\end{enumerate}
}

The point of this proposition is to define a uniform map:
\Def{One says that a map $f$ is a uniform map if for all $\epsilon$ there exists $\eta$ satisfying the conditions of the previous proposition.
}

\Dem{The proof essentially follows the same lines as the proof of \ref{metricmap}:

\begin{enumerate}
\item[$(a) \Rightarrow (b)$] The argument for $(a) \Rightarrow (b)$ in \ref{metricmap} applies in exactly the same way here.

\item[$(b) \Rightarrow (c)$] If $\delta(f_!(U)<\epsilon$, then there exists $q$ such that $\delta( B_q f_!(U) ) < \epsilon$ hence one can take $V = B_q f_!(U)$.

\item[$(c) \Rightarrow (a)$] One has
\[ \Delta_{\eta} = \bigvee_{\delta(U)<\eta} U \times U \]
but for each $U$ such that $\delta(U)<\eta$, there exists $V$ such that $\delta(V)<\epsilon$ and $U \leqslant f^*(V)$, hence

\[ \Delta_{\eta} \leqslant \bigvee_{\delta(V)<\epsilon} f^* V \times f^* V = (f \times f)^* (V \times V) \]
\end{enumerate}

}

}

\block{
\Def{A map between two pre-metric locales is said to be ``compatible with $\triangleleft$" if $U \triangleleft V$ implies $f^* U \triangleleft f^*V$.}

Metric maps and uniform maps are in particular compatible with $\triangleleft$ because if $f$ is uniform and if $\pi_1^* U \wedge \Delta_\epsilon \leqslant \pi_2^*(V)$ then, letting $\eta$ such that

\[ \Delta_\eta \leqslant (f \times f)^* \Delta_\epsilon \]

as we have

\[ (f \times f)^* (\pi_1^*(U ) \wedge \Delta_\epsilon) \leqslant (f \times f)^*\pi_2^* V \]

we obtain

\[ \pi_1^*(f^*U) \wedge \Delta_{\eta} \leqslant \pi_1^*(f^*U)) \wedge (f \times f)^* \Delta_\epsilon \leqslant \pi_2^* f^*V \]
i.e. $f^*U \triangleleft_{\eta} f^* V$
}

\block{\label{isomlemma}\Def{A map $f:X \rightarrow Y$ between two pre-metric locales is called an isometric map if $d(f(x),f(y)) = d(x,y)$, i.e. if $\Delta_q = (f \times f)^*(\Delta_q)$.}

We can easily see (by the same kind of argument that \ref{metricmap}) that this is equivalent to the fact that $\delta(\Lcal) = \delta(f_! \Lcal)$ for all sublocales of $X$.

\Lem{If $f$ is an isometric map $X \rightarrow Y$ then for any locally positive sublocale $\Lcal$ of $X$

\[\Lcal \leqslant f^*(B_q f_! \Lcal ) \leqslant B_q \Lcal \]

}
\Dem{The first inequality immediately follows from the fact that $f_! \Lcal \leqslant B_q f_! \Lcal$.
For the second, as $f_!(\Lcal)$ is locally positive (because of \ref{surjectimpopen}) one can write that

\[ B_q f_! \Lcal = \bigvee_{v \in \Ocal(Y)^{<q} \atop v \wedge f_!(\Lcal) > \emptyset } v. \] 

By \ref{fiberwisevsstrong}, $v \wedge f_!(\Lcal) $ is positive if and only if  $f^*(v) \wedge \Lcal$ is. Also, as $f$ is isometric, for any $v \in \Ocal(Y)^{<q}$ , one has $f^*(v) \in \Ocal(X)^{<q}$. Finally

\[f^*(B_q f_! \Lcal) = \bigvee_{v \in \Ocal(Y)^{<q} \atop f^*(v) \wedge \Lcal >\emptyset } f^*(v) \leqslant \bigvee_{w \in \Ocal(X)^{<q} \atop w \wedge \Lcal > \emptyset} w = B_q \Lcal. \]
}
}

\block{We now consider two toposes $\Ecal$ and $\Tcal$, a geometric morphism $f:\Ecal \rightarrow \Tcal$ and $X$ a pre-metric locale in $\Tcal$.
As $f^\#$ is a functor from locale in $\Tcal$ to locale in $\Ecal$ commuting to projective limit and $f^\#(\Length_{\Tcal}) \simeq \Length_{\Ecal}$, we obtain a map $f^\#(d) : f^\#(X) \times f^\#(X) \rightarrow \Length$. Moreover all the axioms asserting that $d$ is a pre-distance can be pulled back turning $f^\#(X)$ into a pre-metric locale.

\Prop{Let $\Lcal,\Mcal$ be a sublocales of $X$, then (as sublocales of the pre-metric locale $f^\#(X)$) one has:

\begin{itemize}
\item If $\delta(\Lcal)<q$ then $\delta(f^\#(\Lcal))<q$.
\item If $\Lcal \triangleleft_q \Mcal$ then $f^\#(\Lcal) \triangleleft_q f^\#(\Mcal)$.
\item If $\Lcal$ is locally positive then $B_q f^\#(\Lcal) = f^\#(B_q \Lcal)$.
\end{itemize}

}

\Dem{ $f^\#$ is a functor commuting to all projective limits, in particular pull-backs, products and intersections, and by definition of the metric $f^\#(\Delta_q)=\Delta_q$ hence
\[ \Lcal \times \Lcal \subseteq \Delta_{q'} \]
implies \[ f^\#(\Lcal) \times f^\#(\Lcal) \subseteq \Delta_{q'} \]
and
\[ \pi_1^*(\Lcal) \wedge \Delta_q \subseteq \pi_2^*(\Mcal) \]
implies
\[ \pi_1^*(f^\#(\Lcal)) \wedge \Delta_q \subseteq \pi_2^*(f^\#(\Mcal)) \]
which proves the first two points.

The third point is harder because in general the pull-back $f^\#$ does not commute with the direct image functor $(\pi_2)_!$. But if we assume that $\Lcal$ is locally positive, then the map

\[ \pi_1^*(\Lcal) \wedge \Delta_q \rightarrow B_q \Lcal \]

is the restriction of the projection from $\Lcal \times X$ to $X$ and hence is an open map. In particular (as we know that it is a surjection by definition) it is an open surjection and hence its pull-back by $f^\#$ is again an open surjection. In particular, the maps

\[ \pi_1^*(f^\#(\Lcal)) \wedge \Delta_q \rightarrow f^\#(B_q \Lcal) \rightarrow f^\#(X) \]
form a factorisation surjection/inclusion and, by uniqueness of such a factorisation, we obtain the third point.
}

}

\block{\label{Descentpremetric}We also note that if we define $\Ccal(\Tcal)$ to be the category of pre-metric locales and metric maps internal to  $\Tcal$, then open surjections are descent morphisms for $\Ccal$ (see \ref{Descenttheory}) : If $f : \Ecal \rightarrow \Tcal$ is an open surjection and $(X,d)$ is a pre-metric locale in $\Ecal$ endowed with a descent data then it is in particular a descent data on $X$ as a locale, so as locale descend along open surjections, $X$ comes from a locale $X'$ in $\Tcal$. As the $\epsilon : \pi_1^{*} X \rightarrow \pi_2^* X$ is an isomorphism in the category of metric maps it is an isometric map and hence the distance is a morphism in $Des(f,\Ccal)$ and hence also descends into a function $d': X' \times X' \rightarrow \Length$. All the axioms defining a pre-distance are equality relations (and inequality for the specialisation order), hence as they are satisfied by the pull-back of $(X',d')$ along an open surjection they are also satisfied by $(X',d')$. Hence $(X,d)$ is the pull-back of the pre-metric locale $(X',d')$. This proves that the functor $\Tcal \rightarrow Des(f,\Ccal)$ is essentially surjective, but it is also fully faithful for similar reasons: a metric map commuting to descent data is in particular a map of locales commuting to descent data, and as $f$ is an open surjection a map $h$ is metric if and only if $f^*(h)$ is metric. }

\subsection{Metric locales}

\block{If $(X,d$) is a pre-metric locale, then the various properties given in \ref{premetricprop} show that, essentially, the ``topology defined by $d$" (whatever the precise meaning of this is) is coarser than the topology of $X$, but nothing forces them to agree. For example, a metric set in the usual sense (with a distance function taking value in $\Length$), gives a pre-distance on a discrete locale, and the topology defined by $d$ can disagree with the discrete topology. That is why we require the following additional property:

\Def{A Metric locale is a pre-metric locale $X$ such that for all $U \in \Ocal(X)$, 
\[ U = \bigvee_{V \in \Ocal(X) \atop V \triangleleft U} V . \] }

This definition is equivalent to the fact that the family $(B_q V)_{V \in \Ocal(X), q \in \Q_+^*}$ forms a basis of the topology. Indeed $V \triangleleft_q U$ is equivalent to $B_q V \leqslant U$ and $B_q V = \bigvee B_{q'}V$ for $q'<q$, hence this asserts that the open balls form a basis of the topology.

Also if $X$ is metric and $f$ is a geometric morphism then $f^\#(X)$ is also metric because the $B_q V$ for $V\in f^*(\Ocal(X))$ form a basis of $f^\#(X)$.

\Prop{A Metric locale satisfies the following separation axiom: the diagonal embedding
 \[ X \rightarrow \bigwedge_q \Delta_q \] 
is an isomorphism (where the intersection is an intersection of sublocale). 
}

The intuitive reason for this is that if we consider two points $(x,y)$ in $\bigwedge_q \Delta_q $ then by definition $d(x,y)=0$. If the open balls form a basis of the topology then for any open $U$, $x\in U$ if and only if $y \in U$, but for points of a locale this implies that $x=y$. The following proof is just the translation of this argument in terms of generalized points.

\Dem{Consider $f : Y\rightarrow \bigwedge_q \Delta_q$ a map, and let $f_1$ and $f_2$ be the two components $Y \rightarrow X$ of $f$. Let $U$,$V$ be two open sublocales of $X$ such that $U \triangleleft_q V$. Then

\[ \pi_1^*(U) \wedge \Delta_q \leqslant \pi_2^*(V). \]

Applying $f^*$ to each side gives

\[ f_1^*(U) \wedge f^*(\Delta_q) \leqslant f_2^*(V),\]

and as $f^*(\Delta_q)=Y$ by hypothesis, one has $f_1^*(U) \leqslant f_2^*(V)$.

Finally, writing $V = \bigvee_{U \triangleleft V } U$ one has:

\[ f_1^*(V) = \bigvee_{U \triangleleft V } f_1^*(U) \leqslant f_2^*(V).\]
The converse inequality follows by symmetry and hence $f_1=f_2$ i.e. $f$ factors into the diagonal embedding, and this concludes the proof.
}

In particular, as by \ref{overlinedelta}, 

\[ \bigwedge \Delta_q = \bigwedge \overline{\Delta_q} \]

The diagonal embedding of a metric locale is fiberwise closed, one says that metric locales are \emph{fiberwise separated}.

}

\block{\label{isometricsublocale}\Prop{Let $X$ be a metric locale, and $Y$ a pre-metric locale. Let $f : X \rightarrow Y$ be an isometric map. Then $X$ is a sublocale of $Y$ i.e. $f^*$ is onto.

More generally, if we only assume that $X$ is pre-metric then we obtain the inequalities

\[ \forall U \in \Ocal(X), \bigvee_{V \triangleleft U} V \leqslant f^*f_*(U) \leqslant U \]
}

The proposition follows from Lemma \ref{isomlemma}:
\Dem{
Let $U$ be any open sublocale of $X$, such that

\[ U = \bigvee_{V \triangleleft U} V \]

For any $V \triangleleft_q U$ one has by Lemma \ref{isomlemma}

\[V \leqslant f^*(B_q f_! V) \leqslant U \]

hence
\[ U = \bigvee_{q, V \triangleleft_q U} f^*(B_q f_! V) = f^* \left( \bigvee_{q, V \triangleleft_q U } B_q f_! V \right) \]

In particular, if $X$ is metric, then this works for an arbitrary $U$ and $f^*$ is surjective.

If $X$ is no longer metric, then let $U' = \bigvee_{V \triangleleft U} V$, then $U'$ satisfy $U' =\bigvee_{V \triangleleft U'} V$ and hence the first part can be applied to $U'$ and there exists $V$ such that $U' =f^*(V)$. In particular, as $f^*(V) \leqslant U$ we obtain that $V \leqslant f_*(U)$ and hence

\[ U' =f^*(V) \leqslant f^*(f_*(U)). \]

The inequality $f^*(f_*(U)) \leqslant U$ being always true this concludes the proof.

}
}

\block{\label{strongdensity} The following proposition allows one to extend by density relations between continuous functions with values in metric locale.
\Prop{Let $f,g : X \rightrightarrows Y$ be two maps of locales with $Y$ a metric locale (or more generally a fiberwise separated locale). Assume that $f$ and $g$ coincide on some fiberwise dense sublocale $T \subset X$. Then $f=g$. }

\Dem{Let $V$ be the pull-back of the diagonal of $Y$ by the map $(f,g) : X \rightarrow Y \times Y$. As fiberwise closeness is stable under pull-back (see \cite{sketches} C1.2.14(v) ), $V$ is a fiberwise closed sublocale of $X$, containing the fiberwise dense sublocale $T$, hence $V=X$, and this concludes the proof.}

}

\block{\label{extensionofmetricprop} We will also sometimes need to extend by continuity ``metric relations" between functions, which will generally be about comparing functions with value in $\Length$. As $\Length$ is not fiberwise separated, it is not possible to apply directly the previous result. However, one has the following statement:

\bigskip

We will say that a function from $m: X  \rightarrow \Length$ is admissible if there exist two families of functions $f_1, \dots f_n$ and $g_1,\dots,g_n$ from $X$ to pre-metric locales $X_1, \dots X_n$ and a commutative diagram:

\[
\begin{tikzcd}[ampersand replacement=\&]
(\overline{\R_+})^n \arrow{r} \arrow{d} \& \overline{\R_+} \arrow{d} \\
\left( \Length \right)^n \arrow{r}{\lambda} \& \Length \\
\end{tikzcd}
\]

(where the vertical arrows are the canonical maps) such that: 
\[ m(x) = \lambda(d(f_1(x),g_1(x)), \dots, d(f_n(x),g_n(x))).\]

It is probably possible to use a more general definition of ``admissible" map, but this one will be enough for all the applications appearing here.

\Prop{Assume that one has two admissible maps $m_1,m_2: X \rightrightarrows \Length$ such that one has an inequality $m_1 \leqslant m_2$ on some fiberwise dense sublocale $S$ of $X$ a locally positive locale, then the inequality holds one the whole $X$.}

\Dem{The idea is to pull-back everything to some boolean locale $\Bcal$. In the logic of $\Bcal$, thanks to \ref{LEMimpdistcont} the admissible functions $m_1$ and $m_2$ will factor as functions $X \rightrightarrows \overline{\R}$ still satisfying an inequality over $S$. The pull-back of $S$ is still fiberwise dense in the pull-back of $X$ because of \ref{pullbackfwdense}, but, contrary to $\Length$, $\R$ is (fiberwise) separated and hence one can conclude that in the category of sheaves over $\Bcal$ the pull-backs of $m_1$ and $m_2$ agree on the pull-back of $X$ by \ref{strongdensity}. This implies that (in the base topos) one has a diagram:

\[
\begin{tikzcd}[ampersand replacement=\&]
\Bcal \times X \arrow[yshift=0.7ex]{r}{m_1 \leqslant m_2} \arrow[yshift=-0.7ex]{r} \arrow{d}{\pi_2} \& \Bcal \times \Length \arrow{d}{\pi_2} \\
X \arrow[yshift=0.7ex]{r}{m_1,m_1} \arrow[yshift=-0.7ex]{r} \& \Length \\
\end{tikzcd}
\]

In order to conclude that $m_1 \leqslant m_2$ it is enough to choose $\Bcal$ such that $\pi_2 : \Bcal \times X \rightarrow X$ is surjective. It is possible, indeed, if one chooses a boolean locale $\Bcal$ which covers $X$, i.e. with a surjective map $s: \Bcal \twoheadrightarrow X$ then:

\[
\begin{tikzcd}[ampersand replacement=\&]
\Bcal \times \Bcal \arrow{r} \arrow{d}{\pi_2} \& \Bcal \times X \arrow{d}{\pi_2} \\
\Bcal \arrow{r}{s} \& X \\
\end{tikzcd}
\]

The projection $\pi_2 : \Bcal \times \Bcal \rightarrow \Bcal$ is a surjection because it has a section, the map $s : \Bcal \rightarrow X$ is surjective by hypothesis, hence the diagonal map is surjective. This implies that the map $\pi_2 : \Bcal \times X \rightarrow X$ is surjective and hence it concludes the proof.
}

Of course the same result where the inequality is replaced by an equality also holds by two applications of this result.

}

\subsection{Completion of a metric locale}

\blockn{In this subsection we will define the completion of pre-metric locale as the space of minimal Cauchy filters. The same idea has been previously used by S.Vickers in \cite{vickers2005localic}. }

\block{\Def{Let $X$ be a pre-metric locale. A basis $B$ of $X$ is said to be a metric basis if and only if $B$ contains only positive elements, and if $V\in B$ implies $B_qV \in B$.}

This definition can easily be changed without altering the main result of this article, we have chosen it only because it is the simplest notion we have found which is strong enough to assert that the basis will be well behaved and weak enough so that the natural examples we will encounter in practice satisfy this definition, like for example the basis of all open balls on a normed space.

Of course if $B$ is an arbitrary basis of $X$ (composed of positive elements) one can consider the metric basis generated by $B$ by adding to $B$ all the elements of the form $B_{q_1} \dots B_{q_n} V$ for $V \in B$ and $(q_i)$ a finite sequence of positive rational numbers.
Also, if $B$ is a metric basis on $X$ in a topos, then the pull-back of $B$ by any geometric morphism $f :\Ecal \rightarrow \Tcal$ is a metric basis of the pull-back of $X$.
}

\block{\label{cauchyfilterdef}\Def{Let $X$ be a pre-metric locale endowed with a metric basis $B$, a $B$-Cauchy filter on $X$ is a subset $\Fcal \subseteq B$ such that: 

\begin{itemize}
\item[(CF1)]For all $V \in \Fcal$ and $U \in B$ such that $V \leqslant U$ one has $U \in \Fcal$. .
\item[(CF2)]If $U,V\in \Fcal$ then there exists $W \in B$ such that $W\leqslant U$ and $W \leqslant V$ and $W\in \Fcal$.
\item[(CF3)]For all positive rational numbers $q$, there exists $U\in \Fcal$ such that $\delta(U)<q$.
\end{itemize}

A $B$-Cauchy filter is said to be regular if it satisfies additionally:
\begin{itemize}
\item[(CF4)]For all $U\in \Fcal$ there exists $V\in \Fcal$ such that $V \triangleleft U$.
\end{itemize}

A Cauchy filter on $X$ (without specifying the basis) is a $B$-Cauchy filter on $X$, for $B=\Ocal(X)^+$.
}

We insist on the fact that $B$ (as a metric basis) is always assumed to be a subset of $\Ocal(X)^+$. This is why there is no axiom asserting that $\emptyset$ is not an element of $\Fcal$, or that all the elements of $\Fcal$ are positive.
}

\block{\label{regularisationofauchyfilter}\Prop{Any $B$-Cauchy filter $\Fcal$ contains a unique regular Cauchy filter which is $\Fcal^r = \{V \in B | \exists u \in \Fcal, u \triangleleft V \}$.}

\Dem{One easily checks that $\Fcal^r$ is a regular $B$-Cauchy filter. Conversely, let $\Fcal'$ be a regular $B$-Cauchy filter included in $\Fcal$, then for any $U \in \Fcal'$ there exists by $(CF4)$ an element $V \in \Fcal$ such that $V \leqslant B_q V \leqslant U$, hence $U \in \Fcal^r$, which proves that $\Fcal' \subset \Fcal^r$. Let now $U \in \Fcal^r$, by definition there exists $V \in \Fcal$ such that $V \leqslant B_q V \leqslant U$, by $(CF3)$ there exists $W \in \Fcal'$ such that $\delta(W)<q$ and by $(CF2)$ there must be an element $\tau$ of $\Fcal$ such that $\tau \leqslant W$ and $\tau \leqslant V$. In particular, $W \wedge V> \emptyset$ and hence (by the point \ref{PMP_BqLunion} of \ref{premetricprop}) $W \leqslant B_{q} V \leqslant U$ and $U \in \Fcal'$ which concludes the proof. 
}

Hence regular Cauchy filters correspond to the notion of minimal Cauchy filter, this explains why we will later construct the completion of a locale as the classifying space of regular Cauchy filters, by analogy with the classical construction of the completion of a uniform space as a uniform structure on the set of minimal Cauchy filters (see \cite[Chap. II.7]{bourbaki1966elements}).
}

\block{\label{Cauchyfiltermainlemma}\Lem{Let $X$ be a pre-metric locale endowed with a metric basis $B$, and let $\Fcal$ be a regular Cauchy filter on $X$. Then for any $U \in \Fcal$, there exists $V \in B \wedge \Fcal$ such that $V \leqslant U$.}

\Dem{Let $U \in \Fcal$, by $(CF4)$ there exists $U' \triangleleft_q U$ such that $U' \in \Fcal$. Also by $(CF3)$ there exists an element $W \in \Fcal$ such that $\delta(W)<(q/3)$ and as $B$ is a basis and $W$ is positive there exists $b \leqslant W$ with $b \in B$.
Let $V =B_{q/3} b$, then, by the point \ref{PMP_diameterB} of \ref{premetricprop}, one has $\delta(V)< q$ , also $V \in B$ because $B$ is metric, $W \leqslant V$ because $b \wedge W =b$ is positive and $\delta(W)<q/3$ and hence $V \in \Fcal$. Also by $(CF2)$ there exists $V' \in \Fcal$ such that $V' \leqslant V \wedge U'$, as $V'$ is positive this implies that $V \leqslant B_q U' \leqslant U$. As $V \in B \wedge \Fcal$, this concludes the proof.
}
}
\block{\label{CauchyfiltermainlemmaCor}
\Cor{The map $\Fcal \rightarrow B \wedge \Fcal$ induces a bijection between the set of regular Cauchy filters on $X$ and the set of regular $B$-Cauchy filters on $X$.}

We also mention that, as the following proof will show, this proposition holds for any family $B$ satisfying the conclusion of the previous lemma (\ref{Cauchyfiltermainlemma}) even if it is not a metric basis or even if it is not a basis at all.

\Dem{Let $\Fcal$ be a regular Cauchy filter on $X$. We will first prove that $\Fcal' =\Fcal \wedge B$ is a regular $B$-Cauchy filter, this is essentially immediate by Lemma \ref{Cauchyfiltermainlemma}:

\begin{itemize}
\item If $U\leqslant V$ with $V\in \Fcal'$ and $U \in B$ then $U \in \Fcal$ and hence $U \in \Fcal'$ because $\Fcal$ satisfy $(CF1)$.
\item If $U,V \in \Fcal'$ then there exists $W \in \Fcal$ such that $W \leqslant U \wedge V$ and by the lemma there exists $W' \in \Fcal'$ such that $W' \leqslant W \leqslant U,V$.
\item There exists $U \in \Fcal$ such that $\delta(U)<q$ and (by the lemma) a $U' \leqslant U$ such that $U' \in \Fcal'$, hence $\delta(U')<q$.
\item Let $U \in \Fcal'$, there exists $V \in \Fcal$ such that $V \triangleleft U$, then any $V' \leqslant V$ with $V' \in \Fcal'$ (again given by the lemma) works.
\end{itemize}

Now $\Fcal$ can be reconstructed from $\Fcal'$ by the lemma together with $(CF1)$ : 

\[ \Fcal = \{ U | \exists U' \in \Fcal', U' \leqslant U \}. \]

And if you take $\Fcal'$ to be any regular $B$-Cauchy filter, then the previous formula defines a $\Fcal \subseteq \Ocal(X)^+$ which is easily checked to be a regular Cauchy filter as well, and by $(CF1)$ $\Fcal' =\Fcal \wedge B$. This concludes the proof.
}

}

\block{\label{completiondef}Let $X$ be a pre-metric locale, and $B$ be a metric basis on $X$, the theory of regular $B$-Cauchy filters as defined in \ref{cauchyfilterdef} is clearly a propositional geometric theory with basic propositions indexed by $B$. Hence it has a classifying space $\widetilde{X}_B$.

If $X$ is a pre-metric locale in a topos $\Tcal$ and if $f : \Ecal \rightarrow \Tcal$ is a geometric morphism, then $f^\#(\tilde{X}_B) \simeq \widetilde{f^\#(X)_{f^*(B)}}$ because the pull-back of a classifying locale classifies the pull-back of the theory and the pull-back of the theory of regular $B$-Cauchy filter is exactly the theory of regular $f^*(B)$-Cauchy filter on $f^\#(X)$.
But by \ref{CauchyfiltermainlemmaCor} the points of $\widetilde{X}_B$ do not depend on $B$, and hence by the observations we just made, their points on any topos over the base topos do not depend on $B$, and all the $\widetilde{X}_B$ are isomorphic.

\Def{The completion $\widetilde{X}$ of $X$ is defined as the classifying locale $\widetilde{X}_B$ of the theory of regular $B$-Cauchy filters on $X$ for any metric basis $B$ of $X$.}

Also if $U$ is any positive open sublocale of $X$ we denote by $U^\sim$ the open sublocale of $\widetilde{X}$ corresponding to the proposition $``U \in \Fcal"$. It is a general fact about classifying spaces that the $U^\sim$ form a pre-basis of the topology of $X$, but the axiom $(CF2)$ show that for any metric basis $B$ of $X$, the $U^\sim$ with $U \in B$ form a basis of $\widetilde{X}$. If $U$ is not necessarily positive, one can still defined $U^\sim$ by

\[ U^\sim = \bigvee_{V \leqslant U \atop  V > \emptyset } V^\sim. \]

When $U>\emptyset$, the two possible definitions of $U^\sim$ are compatible because

\[ \bigvee_{V \leqslant U \atop  V > \emptyset } V^\sim = U^\sim \]

}

\block{\label{mapstocompletion}\Prop{Let $Y$ be a locale, a morphism $f$ from $Y$ to $\widetilde{X}$ corresponds to a map $\tau : B \rightarrow \Ocal(Y)$ such that:

\begin{enumerate}
\item $\tau$ is non-decreasing.
\item $\displaystyle \tau(U) \wedge \tau(V) \leqslant \bigvee_{W \in B \atop W \leqslant U \wedge V}\tau(W)$
\item $\displaystyle \bigvee_{U \in B \atop \delta(U)<q}(\tau(U))=Y$
\item $ \displaystyle \tau(U) \leqslant \bigvee_{V \in B \atop V \triangleleft U} \tau(V) $
\end{enumerate}
Moreover this correspondence is characterized by the relation $\tau(U) =f^*(U^\sim)$.
Also if $\tau$ only satisfies the first three properties, then there exists a unique $\tau^r$ such that $\tau^r$ satisfy the four properties and $\tau^r \leqslant \tau$ for the pointwise ordering and one has

\[ \tau^r(U) = \bigvee_{V \in B \atop V \triangleleft U} \tau(V) \]

}

\Dem{A morphism from $Y$ to $\widetilde{X}$ is the data of a regular Cauchy filter on $X$ in the internal logic of $Y$. i.e. for each $U \in B$ one should have a proposition $\tau(U) := `` U \in \Fcal "$ satisfying (internally) the axiom $(CF1-5)$. The four properties given for $\tau$ corresponds exactly to the externalisation of the four axioms $(CF1-4)$ (in the right order).

\smallskip

If $\tau$ only satisfies the first three properties then it is just a $B$-Cauchy filter on $X$ and in this case one can apply \ref{regularisationofauchyfilter} and there is a unique regular $B$-Cauchy filter $\tau^r \leqslant \tau$ and it is indeed given by 

\[ \tau^r(U) = \bigvee_{V \in B \atop V \triangleleft U} \tau(V) \]

which is the direct translation of $U \in \tau^r$ if there exists $V \triangleleft U$ with $V \in \tau$.

}

Of course, the inequalities in the axioms $2.$ and $4.$ are in fact equalities because the axiom $1.$ implies the reverse inequalities.

}

\block{\label{injecttocompl}\Prop{There is a map $i$ from $X$ to $\widetilde{X}$ defined by

 \[ i^*(U^{\sim})= \bigvee_{V \triangleleft U} V. \]

Moreover, for any $U \in \Ocal(X)$,

\[ i_*(U) = U^\sim \]

}

\Dem{The inclusion map $e: \Ocal(X)^+ \rightarrow \Ocal(X)$ clearly satisfies the first three points of \ref{mapstocompletion}. Hence the map

\[ e^r(U) = \bigvee_{V \triangleleft U} V \]

satisfies the four points of \ref{mapstocompletion} and hence there is a map $i : X \rightarrow \widetilde{X}$ such that for any $U \in \Ocal(X)^+$ one has $i^*(U^\sim)=e^r(U)$. But as $U^\sim$ is defined as $\bigvee_{V \leqslant U \atop V > \emptyset} V^\sim$ this formula immediately extends to an arbitrary $U$.

We still have to prove that $i_*(U)=U^\sim$. As $i^*(U^\sim) \leqslant U$, one has $U^\sim \leqslant i_*(U)$. Let $V$ an arbitrary open sublocale of $X$ such that $V^\sim \leqslant i_* U$ hence,

\[ \bigvee_{V' \triangleleft V} V' \leqslant U \]

Consider an arbitrary Cauchy filter $F$ on $X$ such that $V \in F$. Then there exists $V' \triangleleft V$ such that $V' \in F$ and hence $U \in F$. This proves that $V^\sim \leqslant U^\sim$ and hence, as $V^\sim \leqslant U^\sim$ imply $V^\sim \leqslant i_*(U)$ one has $V^\sim \leqslant i_* U$ if and only if $V^\sim \leqslant U^\sim$ hence as the $V^\sim$ form a basis of $\widetilde{X}$ this proves that $i_*(U)=U^\sim$.

}
}

\block{\label{XfiberwisedensetildeX}\Prop{The canonical map $i:X \rightarrow \widetilde{X}$ is fiberwise dense and $\widetilde{X}$ is locally positive.}

\Dem{The $(B_q V)^\sim$ for $q$ a positive rational number and $V$ a positive element of $\Ocal(X)$ form a basis of $\widetilde{X}$. Indeed, the $U^\sim$ for $U \in \Ocal(X)^+$ form a basis, and for any $U \in \Ocal(X)$ by $(CF4)$, 

\[ U^\sim = \bigvee_{V \triangleleft U \atop V> \emptyset} V^\sim = \bigvee_{B_q V \leqslant U} (B_q V)^\sim.\]

Moreover, 
\[ i^*((B_q V)^\sim) = \bigvee_{U \triangleleft B_q V} U \geqslant \bigvee_{q' <q} B_{q'} V = B_q V. \]

Hence one has a basis of elements of $\widetilde{X}$ whose pre-image by $i$ are positive. This implies that $\widetilde{X}$ has a basis of positive elements and that for each positive element of $\widetilde{X}$ its pre-image along $i$ is positive, which concludes the proof.
}
}

\block{\label{distcompletion}\Prop{There is a distance function $d$ on $\widetilde{X}$ such that

\[ \Delta_q = \bigvee_{U \in \Ocal(X)^{<q}} U^\sim \times U^\sim. \]
}

One might note that this definition of the distance on $\widetilde{X}$ is the point-free formulation  of the more usual definition:

\[ d(\Fcal,\Fcal')<q\text{ if and only if } \exists u \in \Fcal \wedge \Fcal' \text{ with } \delta(u) <q \]

which is equivalent if interpreted in terms of generalized points.

\Dem{ Let $U \in \Ocal(X)$ such that $\delta(U)<q$. Then there exists $q'$ such that $\delta(U)<q'$ and $U^\sim \times U^\sim \leqslant \Delta_{q'}$. Hence

\[ \Delta_q = \bigvee_{q'<q} \Delta_{q'}, \]

which proves that this formula defines a function $d : \widetilde{X} \times \widetilde{X} \rightarrow \Length$. This function is clearly symmetric, and the diagonal embeddings factor into $\Delta_q$ because the $U^\sim$ with $\delta(U)<q$ cover $\widetilde{X}$ by axiom $(CF3)$. The last point to check is the triangular inequality, but:

\[ \pi_{1,2}^* (\Delta_q) \wedge \pi_{2,3}^*(\Delta_{q'}) = \bigvee_{\delta(U)<q \atop \delta(U')<q'} U^\sim \times (U^\sim\wedge U'^\sim) \times U'^\sim \]

\[ (\pi_{1,3})_!\left(\pi_{1,2}^* (\Delta_q) \wedge \pi_{2,3}^*(\Delta_{q'}) \right) = \bigvee_{\delta(U)<q \atop {\delta(U') < q' \atop U \wedge U' >\emptyset}} U^\sim \times U'^\sim. \]

Since $U^\sim \times U'^\sim \leqslant (U\vee U')^\sim  \times (U\vee U')^\sim $ and as we are restricted to the case $U\wedge U'> \emptyset$, one has $\delta(U \vee U') < q +q'$ by point \ref{PMP_triangulardiamaters} of \ref{premetricprop}, hence  $U^\sim \times U'^\sim  \subset \Delta_{q+q'}$ and

\[(\pi_{1,3})_!\left(\pi_{1,2}^* (\Delta_q) \wedge \pi_{2,3}^*(\Delta_{q'}) \right) \leqslant \Delta_{q+q'}, \]

which is the triangular inequality. The last point to prove is that this pre-distance is a distance. This a consequence of the following lemma.}

\Lem{For any $U \in \Ocal(X)$ one has $B_q (U^\sim) \leqslant (B_q U)^\sim$. In particular, if $U \triangleleft_q V$ then $U^\sim \triangleleft_q V^\sim$.}

\Dem{Indeed, for any $W \in \Ocal(X)$ such that $\delta(W)<q$ and $U^\sim \wedge W^\sim$ is positive, $(CF2)$ proves that $U \wedge W$ is positive, hence, from the definition of $\Delta_q$:

\[ B_q(U^\sim) = (\pi_2)_!(\pi_1^*(U^\sim) \Delta_q) = \left( \bigvee_{\delta(W)<q \atop U^\sim \wedge W^\sim >0} W^\sim \right) \leqslant (B_q U)^\sim  \]
which concludes the proof of the lemma.}

This lemma allows to finish the proof of the proposition, indeed, by $(CF4)$, $V^\sim = \bigvee_{U \triangleleft V} U^\sim$, hence any $V \in \Ocal( \widetilde{X})$ can be written as

\[ V = \bigvee_{U^\sim \leqslant V } U^\sim = \bigvee_{A^\sim \triangleleft U^\sim \leqslant V} A^\sim. \]

}

\block{\label{compextension}\Prop{Let $S \rightarrow Y$ be a fiberwise dense isometric map between two pre-metric locales, let $X$ be any pre-metric locale and $f:S \rightarrow \widetilde{X}$ be a uniform map. Then there exists a unique extension $\widetilde{f} : Y \rightarrow \widetilde{X}$. }

\Dem{The uniqueness of the extension follows from the fact that $\widetilde{X}$ is metric (\ref{distcompletion}) and the result of \ref{strongdensity}, so we only have to prove the existence. We will use \ref{mapstocompletion} for this.
Let $\tau : \Ocal(X)^+ \rightarrow \Ocal(Y)$ defined by:

\[ \tau(U) = i_* f^*(U^\sim) \]
where $i$ denote the embeddings of $S$ into $Y$.

We will first check that $\tau$ satisfies the first three properties of \ref{mapstocompletion}:

\begin{enumerate}

\item $i_*,f^*$ and $U \mapsto U^{\sim}$ are all order preserving. Hence $\tau$ is order preserving.

\item One has $U^\sim \wedge V^\sim = (U \wedge V)^\sim$ (essentially by $(CF2)$) hence as $i_*$ and $f^*$ also commute to binary intersection one has: $\tau(U) \wedge \tau(V) = \tau(U \wedge V)$. This is not enough to conclude immediately the proof of this point because $U \wedge V$ might fail to be positive. Fortunately, if one assumes that $\tau(W)=i_*f^*(W^\sim)$ is positive, then $i^* i_* f^*(W^\sim)$ is also positive because $i$ is fiberwise dense, which implies that $f^*(W^\sim)$ is positive (because it is bigger than $i^* i_* f^*(W^\sim)$) and hence that $W^\sim$ is positive, which finally implies that $W$ is positive (by \ref{XfiberwisedensetildeX} and \ref{injecttocompl}). Hence one can write that

\[\tau(U) \wedge \tau(V) = \tau(U \wedge V)= \bigvee_{\tau(U \wedge V) > \emptyset} \tau(U \wedge V) \leqslant \bigvee_{U \wedge V>\emptyset} \tau(U \wedge V), \]

which proves points $2$.

\item We fix $q$ a positive rational number, and (as $f$ is uniform) $\eta$ such that $\Delta_{\eta} \leqslant (f \times f)^* \Delta_{q/3}$ (see \ref{uniformmap}).

Let $U \in \Ocal(S)^{+,<\eta}$ then (by \ref{uniformmap}) there exists $W \in \Ocal(\widetilde{X})^{<q/3}$ such that $U \leqslant f^*(W)$. 

In particular $W$ is also positive and hence, by $(CF3)$ and the fact that the $V^{\sim}$ form a basis of $\widetilde{X}$, there exists $V_0 \in \Ocal(X)^{+,<q/3}$ such that $V_0^{\sim} \leqslant W$. We define $V = B_{q/3} V_0$. One has $\delta(V) <q$ (by \ref{premetricprop}.\ref{PMP_diameterB}) and $W \leqslant V^\sim$ (by the lemma proved in \ref{distcompletion}), in particular $U \leqslant f^*(V^\sim)$. This proves that

\begin{equation}\label{eq2}
\bigvee_{U \in \Ocal(S)^{+,<\eta}} i_*U \leqslant \bigvee_{V \in \Ocal(X)^{+,<q}}i_*f^*(V^\sim) = \bigvee_{V \in \Ocal(X)^{+,<q}} \tau(V),
\end{equation}

Finally

\[Y = \bigvee_{V \in \Ocal(Y)^{+,<\eta}} V \leqslant \bigvee_{V \in \Ocal(Y)^{+,<\eta}} i_*i^*V = Y. \]
As $i$ is an isometric map, for any $V\in \Ocal(Y)^{<\eta}$ one has $i^*V \in \Ocal(S)^{<\eta}$.
Hence
\begin{equation}\label{eq3}
Y = \bigvee_{V \in \Ocal(Y)^{+,<\eta}} i_*i^*V \leqslant  \bigvee_{U \in \Ocal(S)^{+,<\eta}} i_*U. 
\end{equation}

The inequalities (\ref{eq2}) and (\ref{eq3}) together conclude the proof of the third point.

\end{enumerate}

Hence from \ref{mapstocompletion} there is a map $\tilde{f} : Y \rightarrow \widetilde{X}$ such that $\tilde{f}^*(U^\sim) = \tau^r(U) = \bigvee_{V \triangleleft U} i_*f^* V^{\sim} $. It remains to be proved that $\tilde{f}$ is indeed an extension of $f$, i.e. that $\tilde{f} \circ i = f$.

\[ i^* \tilde{f}^*(U^\sim) = \bigvee_{V \triangleleft U} i^* i_* f^*( V^\sim) \leqslant \bigvee_{V \triangleleft U} f^*(V^\sim) = f^*(U^\sim) \]

Because $\bigvee_{V \triangleleft U} V^\sim = U^\sim$ by $(CF4)$. One the other hand, from the non-metric part of \ref{isometricsublocale}

\[ i^* \tilde{f}^*(U^\sim) = \bigvee_{V \triangleleft U} i^* i_* f^*(V^\sim) \geqslant \bigvee_{V \triangleleft U \atop V' \triangleleft f^*(V^\sim)} V'. \] 

As $f^*$ is uniform it is compatible with $\triangleleft$, hence the set of $V'$ appearing in the last union contains all the $f^*(W^\sim)$ for $W \triangleleft V$ hence

\[ i^* \tilde{f}^*(U^\sim) \geqslant \bigvee_{V \triangleleft U \atop W \triangleleft V } f^*(W^\sim) = f^*(U ^\sim), \]

which proves $i^* \tilde{f}^*(U^\sim) = f^*(U^\sim)$ and concludes the proof.

}

We also note that if the map $f$ is metric (resp. isometric), the extension $\widetilde{f}$ will also be metric (resp. isometric) by an application of \ref{extensionofmetricprop}.

}

\block{\label{compdef}\Th{Let $X$ be a pre-metric locale, then the following conditions are equivalent: 
\begin{enumerate}
\item The map $X \rightarrow \widetilde{X}$ is an isomorphism;
\item $X \simeq \widetilde{Y}$ for some $Y$;
\item For any $S \rightarrow Y$ a strongly dense isometric map between pre-metric locales, and any map from $S$ to $X$ there exists a map from $Y$ to $X$ making the triangle commute;
\item Any strongly dense isometric map from $X$ to a metric locale $Y$ is an isomorphism.
\end{enumerate}

A locale satisfying these conditions is called a {\it complete} metric locale.

}
\Dem{$1. \Rightarrow 2.$ is clear. 

$2. \Rightarrow 3.$ is a direct consequence of \ref{compextension}.

$4. \Rightarrow 1.$ is also clear because the map from $X$ to $\widetilde{X}$ is a dense isometric map.

$3. \Rightarrow 4.$ remains to be proved.
Let $f : X \rightarrow Y$ be a strongly dense isometric map. The identity map from $X$ to $X$ can be extended into a map $g$ from $Y$ to $X$ by $3.$, such that $g \circ f = Id_X$. As, $f \circ g$ restricted to $X$ is the inclusion from $X$ to $Y$, $f \circ g$ is the identity of $Y$ by fiberwise density of $X$ into $Y$ and fiberwise separation of $Y$ ( \ref{strongdensity}) hence $g$ is an inverse for $f$, and they are isomorphisms.

}

It is immediate from point $3.$ that a locally positive fiberwise closed sublocale of a complete locale is also complete.
}

\block{\label{Compltnessdescend}\Prop{If $X$ is a pre-metric locale in a topos $\Tcal$ and $f:\Ecal \rightarrow \Tcal$ is an open (or proper) surjection such that $f^\#(X)$ is complete then $X$ is complete. 
}

\Dem{The pull-back along $f$ of the canonical map $X \rightarrow \widetilde{X}$ is the canonical map $f^\#(X) \rightarrow \widetilde{f^\#(X)}$. Hence as $f^\#$ is a descent functor for the categories of locales, it is in particular conservative and if the pull-back map is an isomorphism, the map $X \rightarrow \widetilde{X}$ is also an isomorphism.}

An immediate corollary of this result is that if $\Ccal(\Tcal)$ is the category of complete metric locales and metric maps between them then objects of $\Ccal$ descend along open surjections. Indeed, it is a full subcategory of the category of pre-metric locales, for which open surjections are descent morphisms as observed in \ref{Descentpremetric}, and this just states that $(X',d')$ is complete if it descends from a complete locale $(X,d)$.

}

\block{\Prop{Let $X$ be a pre-metric locale and let $X_d$ be the regular image of $X$ into $\widetilde{X}$ then $\Ocal(X_d)$ identifies with the set of $U \in \Ocal(X)$ such that

\[U= \bigvee_{V \triangleleft U} V\]

and any map compatible with $\triangleleft$ from $X$ to a metric locale $Y$ factors into $X_d$.
}

\Dem{The regular image of $i:X \rightarrow \widetilde{X}$ is identified as a frame with the image of $i^*: \Ocal(\widetilde{X}) \rightarrow \Ocal(X)$ which is clearly (by \ref{injecttocompl}) the set of open sublocales defined in the proposition. If one has any map $f$ from $X$ to a metric locale $Y$ compatible with $\triangleleft$ then for any $U \in  \Ocal(Y)$,

\[ U = \bigvee_{V \triangleleft U} V\]

Hence, 

\[ f^*(U) = \bigvee_{V \triangleleft U } f(V)^* \]

as $f^*(V)\leqslant f^*(U)$ this proves that $f^*(U) \in \Ocal(X_d)$. Hence $f$ factors into $X_d$.
}

} 

\subsection{Product of metric locales}
\label{metricproduct}
\block{Let $\Lcal$ and $\Mcal$ be two pre-metric locales, one defines a pre-distance on $\Lcal \times \Mcal$ in the following way: $\Delta_q^{\Lcal \times \Mcal} \subset (\Lcal \times \Mcal) \times (\Lcal \times \Mcal)$ is the intersection of the pull-back $\pi_{1,3}^*(\Delta_q^{\Lcal})$ and $\pi_{2,4}^*(\Delta_q^{\Mcal})$ (where the exponent on $\Delta$ indicate to which locale it is related). This corresponds to taking $d((l,m),(l',m')) = \max (d(l,l'),d(m,m'))$, and the classical argument can be adapted (in terms of generalised points) to prove that this is indeed a pre-distance on $\Lcal \times \Mcal$.

\Prop{ $\Mcal \times \Lcal$ endowed with the previously constructed distance function is the categorical product of $\Mcal$ and $\Lcal$ in the category of pre-metric locales and metric maps. }

\Dem{The projection $\pi_1 : \Lcal \times \Mcal \rightarrow \Lcal$ satisfies $\Delta_q \subset \pi_1^*( \Delta_q)$ by construction of the distance function on $\Lcal \times \Mcal$, hence it is a metric map. In particular if $f: X \rightarrow \Mcal \times \Lcal$ is a metric map then the two component $f_1$ and $f_2$ are metric maps. Conversely, assume that $f_1$ and $f_2$ are metric maps. Then

\[ (f \times f)^* (\Delta_q^{\Lcal \times \Mcal}) = (f \times f)^* ( \pi_{1,3}^*(\Delta_q^{\Lcal}) \wedge \pi_{2,4}^* (\Delta_q^{\Mcal})). \]

But $\pi_{1,3}( f \times f) = f_1 \times f_1$ and $\pi_{2,4} (f \times f) = f_2 \times f_2$, hence,

\[ (f \times f)^* (\Delta_q^{\Lcal \times \Mcal}) = (f_1 \times f_1)^*( \Delta_q^{\Lcal}) \wedge (f_2 \times f_2)^*(\Delta_q^{\Mcal}) \]

As we assume that both $f_1$ and $f_2$ are metric,

\[ \Delta_q^X \subset (f_1 \times f_1)^*( \Delta_q^{\Lcal}) \wedge (f_2 \times f_2)^*(\Delta_q^{\Mcal}), \]

This proves that $f$ is also metric and concludes the proof of the proposition.

}

}

\block{\Prop{The product of two complete metric locales is a complete metric locale. More generally the completion of $\Lcal \times \Mcal$ is canonically isomorphic to $\widetilde{\Lcal} \times \widetilde{\Mcal}$. }

\Dem{ Assume that $\Lcal$ and $\Mcal$ are complete. Let $S \rightarrow Y$ be a strongly dense map, and let $f :S \rightarrow \Lcal \times \Mcal$ be an isometric map. Then by the previous result and Proposition \ref{compextension} there is a map $\widetilde{f}: Y \rightarrow \Lcal \times \Mcal$ extending $f$. Hence $\Lcal \times \Mcal$ is complete.

For the second part, $\Lcal \times \Mcal \rightarrow \widetilde{\Lcal} \times \widetilde{\Mcal}$ is a fiberwise dense isometric map with $\widetilde{\Lcal} \times \widetilde{\Mcal}$ complete, hence $\widetilde{\Lcal} \times \widetilde{\Mcal}$ is the completion of $\Lcal \times \Mcal$.
}
}

\subsection{The locale $[X,Y]_1$ of metric maps}
\label{mapingspace}
\blockn{In this subsection we show that it is possible to construct a classifying space $[X,Y]_1$ of metric maps between two metric locales $X$ and $Y$, at least when $Y$ is complete. The key observation underlying this construction is that (in a classical settings) on the set of metric functions the topology of point-wise convergence on any dense subsets is equivalent to the compact-open topology, and that when we endow this set of metric functions with this topology the composition law is bi-continuous. This suggests that this topology classifies metric functions. The general idea of this section is to give a point-free formulation of this topology, by replacing the basic open $``f(x) \in V"$ by $``U \wedge f^{-1}(V)> \emptyset"$ for $U$ a small neighborhood of $x$. }

\block{\Def{Let $X$ and $Y$ be two pre-metric locales. Let $A$ be a basis\footnote{One can actually see that we do not even need $A$ to be a basis. All we need is that for all positive rational $q$ the set of $a \in A$ such that $\delta(a)<q$ cover $X$.} of positive open of $X$ and $B$ be a metric basis of $Y$. We define $[X_A,Y_B]_1$ as the classifying space of the propositional geometric theory on propositions $(U,V)$ for $U \in A$ and $V \in B$ with the axioms:

\begin{itemize}

\item[(MM1)] For all $U' \leqslant U$ and $V' \leqslant V$
\[ (U',V') \vdash (U,V) \]
\item[(MM2)] For all $V \in B$,$U \in A$ and any positive rational number $q$ one has

\[ (U,V) \vdash \bigvee_{u \leqslant U \atop \delta(u)<q} (u,V); \]

\item[(MM3)] For all $U \in A$ and all $q$ positive:

\[ \vdash \bigvee_{V \in B \atop \delta(V)<q} (U,V);\]

\item[(MM4)]For all $U \in A$, $V \in B$ 

\[(U,V) \vdash \bigvee_{V' \in B \atop V' \triangleleft V} (U,V'); \]

\item[(MM5)] Let $W_1,W_2, \tau \in A$, $q_1,q_2 \in \Q$, $V_1,V_2,V_1',V_2' \in B$ such that

\[ \begin{array}{c c}
\delta(W_1)<q_1 & \delta(W_2) <q_2 \\
V_1' \triangleleft_{q_1} V_1 & V_2' \triangleleft_{q_2} V_2 \\
\tau \leqslant W_1 & \tau \leqslant W_2
\end{array} \]

then

\[(W_1,V_1')\wedge (W_2,V_2') \vdash \bigvee_{V \in B \atop V \leqslant V_1 \wedge V_2} (\tau,V) \]

\item[(MM6)]

\[ (U,V) \wedge (U,V') \vdash \delta(V \vee V') \leqslant \delta(U)+\delta(V)+\delta(V'). \]

\end{itemize}

}
}
\block{\label{mappingspacemainth}

The main result of this section is 

\Th{The locale $[X_A,Y_B]_1$ we just constructed does not depend on $A$ and $B$ and classifies metric maps between $X$ and $\widetilde{Y}$. With the propositions $(U,V)$ corresponding to $U \wedge f^*(V^\sim) > \emptyset$. This locale will be denoted $[X,Y]_1$ }

Its proof will occupy us for the rest of this subsection.

}

\block{If $f$ is a geometric morphism from $\Ecal$ to $\Tcal$, then, by the same argument as in \ref{completiondef}:

\[ f^\# ([X_A,Y_B]_1) \simeq [f^\#(X)_{f^*(A)}, f^\#(Y)_{f^*(B')}]_1 \]

So it suffices to show that the points of $[X_A,Y_B]_1$ correspond to metric functions from $X$ to $\widetilde{Y}$ to obtain the announced result.
}

\block{\label{pointfrommap}\Prop{Let $f:X \rightarrow \widetilde{Y}$ be a metric map and let:

 \[ (U,V)_f:= `` U \wedge f^*(V^\sim)>\emptyset " \]

For $U \in A$ and $V \in B$. 
Then this defines a point of $[X_A,Y_B]_1$.
}

\Dem{Axiom $(MM1)$ is immediate. $(MM2)$ holds because for any $V \in B, U \in A$, if $f^*(V^\sim) \wedge U$ is positive then one can write $U$ as a union of $u \in A$ such that $u \leqslant U$ and $\delta(u)<q$ and the locale positivity of $X$ allows one to conclude. Axiom $(MM3)$ and $(MM4)$ hold because the corresponding unions holds in $\widetilde{Y}$.

\bigskip

We now prove axiom $(MM5)$. Let $W_1,W_2,\tau,q_1,q_2,V_1,V'_1,2_2,V'_2$ satisfying the hypothesis of $(MM5)$. We also assume that $(W_1,V'_1)_f$ and $(W_2,V'_2)_f$ holds. Then as $f$ is metric and $V'_i \triangleleft_{q_i} V_i$ then $V_i'^\sim \triangleleft_{q_i} V_i^\sim$ one has

\[ f^*(V_i'^\sim) \triangleleft_{q_i} f^*(V_i^\sim). \]

As $\delta(W_i)<q_i$ and $W_i \wedge f^*(V_i) > \emptyset$ this implies that

\[ W_i \subseteq f^*(V_i^\sim), \]
and hence, as $\tau \leqslant W_1 \wedge W_2$, that 
\[ \tau \subseteq f^*(V_1^\sim \wedge V_2^\sim). \] 

As $\tau$ is positive (the presentation of $X$ is assumed to be locally positive) and $V_1^\sim \wedge V_2^\sim$ is covered by the $V^\sim$ for $V \subseteq V_1 \wedge V_2$ this concludes the proof of $(MM5)$.

\bigskip

We now prove $(MM6)$. Let $U,V$ and $V'$ such that $U \wedge f^*(V^\sim)>\emptyset$ and $U \wedge f^*(V'^\sim)>\emptyset$. Let $q$ and $q'$ such that $\delta(V)<q$ and $\delta(V')<q'$. Let also $\epsilon$ be a positive rational number such that $\delta(V)<q-2\epsilon$ and  $\delta(V')<q'-2\epsilon$. Let $W=B_\epsilon V$ and $W' =B_\epsilon V'$, in particular $\delta(W)<q$ and $\delta(W')<q'$.

One has, by the assumption on $V$ and $V'$ and the fact that $f$ is metric  (see \ref{metricmap} proposition $(c)$):

\[ \delta(W^\sim \vee W'^\sim) \subseteq \delta(W^\sim)+\delta(W'^\sim)+\delta(U) \]

Let $i$ be the isometric map $Y \rightarrow \widetilde{Y}$ of \ref{injecttocompl}, i.e.
\[i^*(V^\sim)=\bigvee_{U \triangleleft V} U. \]

In particular, as $W$ and $W'$ are open balls, one has $i^*(W^\sim)=W$ and $i^*(W'^\sim)=W'$, and $i^*(W^\sim \vee W'^\sim)=W \vee W'$, and as $i$ is isometric, this implies that $\delta(W \vee W') \leqslant \delta(W^\sim \vee W'^\sim)$.

Moreover since $\delta(W)<q$ then by definition of the distance on $\widetilde{Y}$, $W^\sim \times W^\sim \subseteq \Delta_q$, and hence $\delta(W^\sim) \leqslant q$. One deduces from this that

\[ \delta(V \vee V') \leqslant \delta(W \vee W') \leqslant \delta(W^\sim \vee W'^\sim) \leqslant  \delta(W^\sim)+\delta(W'^\sim)+\delta(U)  \leqslant q+q'+\delta(U), \]

which concludes the proof as it has been done for arbitrary $q$ and $q'$ bigger than $\delta(V)$ and $\delta(V')$.

}
}

\block{\Def{To any point $p$ of $[X_A,Y_B]_1$ we associate the function $\tau_p: B \rightarrow \Ocal(X)$ defined by:

\[ \tau_p(V) :=\bigvee_{\delta(W) <q \atop {V' \triangleleft_q V \atop p \in (W,V')}} W  \]

where $V'$ runs through elements of $B$, $W$ through elements of $A$, and $q$ through positive rational numbers.
}

\Prop{If $f$ is a metric map from $X$ to $\widetilde{Y}$ and $p$ is the point of $[X_A, Y_B]$ associated to $f$ in \ref{pointfrommap} then

\[\tau_p(V)=f^*(V^\sim). \]

}
\Dem{ One has by definition:
\[ \tau_p(V) = \bigvee_{\delta(W) <q \atop {V' \triangleleft_q V \atop f^*(V'^\sim) \wedge W > \emptyset }} W. \] 

Hence, as for any $W$ appearing in the supremum one has $W \leqslant f^*(V^\sim)$, we obtain that $\tau_p(V) \leqslant f^*(V^\sim)$.

Conversely,

\[ f^*(V^\sim) = \bigvee_{V' \triangleleft_q V} f^*(V'^\sim) = \left( \bigvee_{V' \triangleleft_q V \atop {\emptyset <W \leqslant f^*(V'^\sim) \atop \delta(W)<q}} W \right) \leqslant \tau_p(V'^\sim). \]

}

}

\block{\Lem{Let $p$ be any point of $[X_A,Y_B]_1$, then:

\[ p \in (U,V) \Leftrightarrow U \wedge \tau_p(V)> \emptyset\] 

}

\Dem{Assume first that $\tau_p(V) \wedge U > \emptyset$. Then there exists $W$ and $V'$ such that $\delta(W)<q$, $V' \triangleleft_q V$ , $(W,V')$ and $W \wedge U>\emptyset$. Applying $(MM5)$, one obtains that there exists $V'' \leqslant V$ such that $p \in (W \wedge U, V'')$ and hence $p \in (U,V)$.

Conversely assume that $p \in (U,V)$, then (by $(MM4)$) there exists $V' \in B$ and a positive $q$ such that $V' \triangleleft_q V$ and $p \in (U,V')$. Also by $(MM2)$ there exists $W \in A$ such that $\delta(W)<q$ and $p \in (W,V')$. But this implies that $W \leqslant  \tau_p(V)$ and as $W \leqslant U$ and $W > \emptyset$ one concludes that $U \wedge \tau_p(V)>\emptyset$.
}

}

\block{At this point, all that remains to be checked in order to prove \ref{mappingspacemainth} is that for any point $p$, $\tau_p$ extends into a map from $X \rightarrow \widetilde{Y}$ and that this map is indeed metric. 

\Prop{The map $\tau_p : B \rightarrow \Ocal(X)$ satisfies the four conditions of \ref{mapstocompletion} and in particular there is a (unique) map $f:X \rightarrow \widetilde{Y}$ such that $f^*(V^\sim)= \tau_p(V)$. }

\Dem{ We recall that

\[ \tau_p(V) :=\bigvee_{\delta(W) <q \atop {V' \triangleleft_q V \atop p \in (W,V')}} W  \]

Also the point $p$ being fixed, we will write $\tau$ instead of $\tau_p$ and $(U,V)$ instead of $p \in (U,V)$.

\begin{enumerate}
\item if $U \leqslant V$ then any $W$ appearing in the supremum defining $\tau(U)$ also appears in the one defining $\tau(V)$ with the same $V'$ and $q$. Hence $\tau$ is order preserving.

\item
\[ \tau(V_1) \wedge \tau(V_2) = \bigvee W_1 \wedge W_2 \]

where the union runs over all $W_1,W_2 \in A$ such that there exist $q'_1,q'_2$ positive rational numbers, and $V'_1 , V'_2 \in B$ such that
\[ \delta(W_i)<q'_i; \]
\[V'_i \triangleleft_{q'_i} V_i; \]
\[ (W_i,V'_i). \]

For any such $W_1$ and $W_2$ there exists a positive rational number $\epsilon$ such that $\delta(W_i) < q'_i - \epsilon$. Let $q_i= q'_i - \epsilon$. One has in particular $\delta(W_i)<q_i$ and

\[V'_i \triangleleft_{q_i} B_{q_i} V'_i \triangleleft_{\epsilon} V_i. \]

Moreover $W_1 \wedge W_2$ can be written as the union of $\tau \in A$ such that $\tau \leqslant W_1 \wedge W_2$ and $\delta(\tau)<\epsilon$. Finally, one can apply $(MM5)$ (taking $B_{q_i} V'_i$ instead of $V_i$) to obtain that there exists $V$ such that

 \[ V \leqslant (B_{q_1} V'_1 \wedge B_{q_2} V'_2) \triangleleft_{\epsilon} V_1 \wedge V_2 \]
and
\[ (\tau,V). \]

This proves that $\tau \leqslant \tau(B_\epsilon V)$ with $B_\epsilon B \leqslant V_1 \wedge V_2$ and $B_\epsilon V \in B$ because $B$ is metric, and hence concludes the proof that.

\[ \tau(V_1) \wedge \tau(V_2) \leqslant \bigvee_{V \in B \atop V \leqslant V_1 \wedge V_2} \tau(V). \]

\item Let $q$ be any positive rational number. Let $W \in A$ such that $\delta(W)<q/3$.
Then by  $(MM3)$ there exists $V' \in B$ such that $\delta(V')<q/3$ and $(W,V')$. Let $V = B_{q/3} V' \in B$, one has: $\delta(W)<q/3$, $ V' \triangleleft_{q/3} V$, $(W,V')$, hence $W \leqslant \tau(V)$ with $\delta(V)<q$ this proves that

\[ W \leqslant \bigvee_{V \in B \atop  \delta(V)<q} \tau(V) \]

As we have done this for an arbitrary $W$ with $\delta(W)<q/3$ this concludes the proof.

\item Let $V \in B$, let $W$ appearing in the union defining $\tau(V)$, i.e. there exists a positive rational $q$, and a $V' \in B$ such that $\delta(W)<q$ and $V' \triangleleft_q V$.

But, there exists a positive rational number $\epsilon$ such that $\delta(W)<q-\epsilon$, and $V' \triangleleft_{q-\epsilon} B_{q-\epsilon} V' \triangleleft_\epsilon V$. Hence

 \[ W \leqslant \tau(B_{q-\epsilon} V' \leqslant \bigvee_{U \in B \atop U \triangleleft V} \tau(U). \]

Finally, we obtain

\[ \tau(V) \leqslant \bigvee_{U \in B \atop U \triangleleft V} \tau(U). \]
\end{enumerate}

}

The fact that the map $f$ induced by $\tau_p$ is metric follows from axiom $(MM6)$ using the characterization $(c)$ of metric maps given in \ref{metricmap}, hence this concludes the proof of theorem \ref{mappingspacemainth}.

}

\subsection{Case of metric sets}
\label{Metricset} 

\block{We define a (pre)metric set as set $X$ endowed with a distance function $d : X \times X \rightarrow \Length$ satisfying the usual axioms for a (pre)distance:

\begin{itemize}
\item $d(x,x) = 0$
\item $d(x,y)= d(y,x)$
\item $d(x,z) \leqslant d(x,y)+d(y,z) $
\end{itemize}

With additionally, $d(x,y)=0 \Rightarrow x=y$ for a metric set.

A (pre)metric set can be seen as a pre-metric locale by seeing its underlying set as a discrete locale. It is in general not a metric locale even if we start with a metric set.
}

\block{We will say that a metric set $(X,d)$ is complete if the natural map $i:X \rightarrow \widetilde{X}$ identifies $X$ with the points of $\widetilde{X}$. As points of $\widetilde{X}$ identify with regular Cauchy filters one easily checks that this is equivalent to the usual (Cauchy filter based) definition of completeness.}

\block{\label{eqweakspatial}\Th{There is an equivalence of categories between the category of weakly spatial complete metric locales (with metric maps) and complete metric sets (with metric maps). }

\Dem{The functors are given by the following construction: to a complete metric set $X$ one associates its localic completion $\widetilde{X}$, which is weakly spatial, because $X$ is fiberwise dense in it, and to a weakly spatial complete metric locale one associates its set of points endowed with the induced distance. These two constructions are functorial on metric maps.

By definition of a complete metric set it identifies with the set of points of its localic completion, and conversely, if $\Lcal$ is a weakly spatial complete metric locale and $X$ is its set of points endowed with the induced distance, then $X \rightarrow \Lcal$ is a fiberwise dense isometric map from $X$ to a complete locale, hence $\Lcal$ is isomorphic to the completion of $X$. This proves that the two functors are inverse from each other on objects. They are also inverse of each other on morphisms, tautologically on one side and by \ref{strongdensity} on the other side.
}
}

\block{\label{netoffunction}The internal application of the fact that the set of points of a complete metric locale is complete in the classical sense can prove directly a result of completeness of the space of functions with values in a complete locale for the uniform distance. This cannot be stated directly in terms of completeness of some metric locale because in general (if the initial space is not locally compact) the space of functions is not a locale, but one has:

\Prop{Let $(f_i)_{i\in I}$ be a Cauchy net of functions between two locales $X$ and $Y$, with $Y$ a complete metric locale. This means that $I$ is a directed (filtering) ordered set and that for all positive rational number $\epsilon$ there exists $i_0 \in I$ such that $\forall i,j \geqslant i_0$, the map $(f_i,f_j)$ factors into $\Delta_{\epsilon} \subset Y \times Y$.

Then the net $f_i$ converges to some (uniquely defined) function $f :X \rightarrow Y$. This mean that there is a unique function $f:X \rightarrow Y$ such that for all positive rational number $\epsilon$ there exists $i_0 \in I$ such that $\forall i \geqslant i_0$, the map $(f,f_i)$ factors into $\Delta_\epsilon$.
 }
 
\Dem{ The net of functions $f_i : X \rightarrow Y$ can be interpreted as a net of points of $p^\# Y$ in the logic of $X$ (where $p$ is the map $X \rightarrow *$). And the fact that it is externally a Cauchy net immediately gives that it is internally a Cauchy net. The usual proof that completeness by filter imply completeness by net is completely constructive\footnote{On the contrary, the converse relies on the axiom of choice.} and hence the fact that $p^\# Y$ is complete implies the convergence of the net $f_i$. Uniqueness of the limit implies that the limit is a global point of $p^\# Y$ in $X$, and hence a map from $X$ to $Y$. One then easily check that the internal convergence together with the external Cauchy condition imply the external convergence. }
}

\block{\label{stackification}In particular the category of complete metric sets identifies with the full subcategory of the category of complete metric locales composed of weakly spatial locales, and by \ref{potentiallyweaklyspatial} any complete metric locale becomes weakly spatial (hence identifies with a complete metric set) after a pull-back to some open locale. We already mentioned that if one defines $\Ccal(\Tcal)$ as the category of complete metric locales over $\Tcal$, then, it is a stack for the topology whose covering are open surjections.

From these observations one can deduce that the stack of internal complete metric locales is the stackification (the analogue of sheafication for stack and pre-stack) of the pre-stack of complete metric sets, that is the universal extension of the notion of complete metrics sets for the descent properties along open surjection.

\bigskip

At this point one could obtain the localic Gelfand duality of \ref{gelfand} directly by observing that the notion of compact regular locale is obtained as the stackification of the notion of compact completely regular locale, and apply the constructive Gelfand duality between compact regular locale and $C^*$ algebra to show that the two pre-stacks are equivalent. This will also avoid the use any of the material of section \ref{mapingspace}, but it will give an extremely uncomfortable definition of the spectrum of a localic $C^*$ algebra. This is why we prefer explicitly constructing the spectrum (in \ref{spectrum}, using the construction of \ref{mapingspace}) before applying the descent argument to show the Gelfand duality.
}

\section{Banach locales and $C^*$ locales}
\label{sectionBanach}
\subsection{Banach locales and completeness}

\block{ \Def{A pre-Banach locale is a locally positive locale $\Hcal$ endowed with:

\begin{itemize}
\item A commutative group law: $+: \Hcal \times \Hcal \rightarrow \Hcal$, with neutral element $0 : * \rightarrow \Hcal$ and an inversion: $x \mapsto -x : \Hcal \rightarrow \Hcal$.
\item An action of $\Q[i]$ (endowed with the discrete topology), $\Q[i] \times \Hcal \rightarrow \Hcal$, satisfying the usual axioms of a (unital) module.
\item A norm function $\Vert . \Vert : \Hcal \rightarrow \Length$
\end{itemize}

where the norm function is expected to satisfy the following conditions:
\begin{itemize}
\item $\forall x,y \in \Hcal \Vert x+y \Vert \leqslant \Vert x \Vert + \Vert y \Vert$
\item $\forall \lambda \in \Q[i], \forall x \in \Hcal, \Vert \lambda x \Vert = |\lambda| \Vert x \Vert$
\item $\Vert 0 \Vert = 0$
\item $\Hcal = \bigvee_{n \in \N} \{ x | \Vert x \Vert < n \} $
\end{itemize}
}

Of course, all the conditions stated in this definition have to be interpreted either in diagrammatic terms or in terms of generalized elements. 

}

\block{\label{prebanachprop}
\Prop{Let $(\Hcal, \Vert .\Vert)$ be a pre-Banach locale. Let $s$ and $p$ denote the maps $\Hcal \times \Hcal \rightarrow \Hcal$ defined by:
\[ s(x,y)=x-y \]
\[ p(x,y)=x+y \]
Let $m$ denote the map $x \mapsto -x$ and $n$ be the norm map, $n : \Hcal \rightarrow \Length$.

Finally we will denote $B_q0 = n^*([0,q[)$ (point \ref{PBP_coherence} ensures that there is no possible confusion).

Then, one has the following facts:
\begin{enumerate}
\item \label{PBP_distdef} The map $n \circ s$ is a pre-distance on $\Hcal$.
\item \label{PBP_openplus}The maps $s$ and $p$ are open maps.
\item The open sublocales $\Delta_q$ coincide with $s^*(B_q 0 )$.
\item If $\Lcal$ is any sublocale of $\Hcal$ then $B_q \Lcal$ coincide with both $p_!(\Lcal \times B_q 0)$ and $s_!(\Lcal \times B_q 0)$.
\item \label{PBP_coherence} $B_q 0$ is the same things as $ B_q \{0 \}$.
\end{enumerate}

}

\Dem{
\begin{enumerate}
\item A proof by generalized points will be exactly the same as the usual proof that $d(x,y)=\Vert x-y \Vert$ is a distance on a normed space.

\item We will consider two maps $\Hcal \times \Hcal \rightarrow \Hcal \times \Hcal$ given by

\[ \tau_p = (p, m \circ \pi_1); \]
\[ \tau_s = (\pi_1,s). \]

These maps correspond in term of generalized points to the maps $\tau_p(x,y)=-x+y,-y)$ and $\tau_s(x,y)=(x,x-y)$, and they are both involutive and hence bijective. The maps $s$ and $p$ are then obtained as $\pi_2 \circ \tau_s$ and $\pi_1 \circ \tau_p$, but as $\Hcal$ is locally positive, both $\pi_1$ and $\pi_2$ are open maps. Hence by composition $s$ and $p$ are open maps.

\item $\Delta_q$ is by definition $d^*([0,q[)$, but as $d= n \circ s$, one has $\Delta_q= s^* n^*([0,q[) = s^*(B_q 0)$.

\item The involutive map $\tau_s$ introduced in the proof of point \ref{PBP_openplus} exchange $\pi^*( \Lcal) \wedge \Delta_q$ with $\Lcal \times B_q0$, indeed:

\[ \tau_s^*( \Lcal \times B_q 0 ) = pi_1^*( \Lcal) \wedge s^*(B_q0) = \pi_1^*(\Lcal) \wedge \Delta_q. \]

Hence $\pi_{2!}(\pi_1^*(\Lcal) \wedge \Delta_q) = (\pi_2 \circ\tau_s)_!( \Lcal \times B_q 0)$ and $\pi \circ \tau_s = s$, which shows that $B_q \Lcal = s_!(\Lcal \times B_q 0)$.

It also coincides with $p_!(\Lcal \times B_q 0)$ because as $n \circ m= n$ one has $m^*(B_q 0)= B_q0$, and as $s = p \circ (Id,m)$ this concludes the proof.

\item From the previous result, $B_q \{ 0 \}$ identifies with $p_!(\{0 \} \times B_q 0)$ but $p$ acts on $\{0 \} \times B_q 0$ as the inclusion of $B_q 0$ in $\Hcal$ (this is the definition of $0$ being the neutral element), hence $p_!(\{0 \} \times B_q 0) = B_q 0$ and this concludes the proof.

\end{enumerate}

}

}

\block{\label{Banachdef}\Prop{Let $\Hcal$ be a pre-Banach locale, the following conditions are equivalent:
\begin{enumerate}
\item[(LB1)] The open sublocales $B_q 0$ form a basis of neighborhoods of $0$.

\item[(LB2)] $\Hcal$ is metric for the distance induced by $\Vert . \Vert$.
\end{enumerate}
}

A pre-Banach locale satisfying either $(LB1)$ or $(LB2)$ is called a Banach locale, we will soon see that there is no need for a completeness assumption: it will be automatic.

\Dem{We will use the same notation $s,p$ as in proposition \ref{prebanachprop}.
Assume $(LB1)$, and let $U$ be any open of $\Hcal$. Consider the open sublocale $p^* U \subset \Hcal \times \Hcal$, and decompose it as a union of basic open sublocales

\[ p^* U = \bigvee_{i \in I} A_i \times B_i \]

where $A_i$ and $B_i$ are open sublocales of $\Hcal$.
Let $i$ such that $(A_i \times B_i) \wedge U \times \{0 \}$ is positive. Then $B_i \wedge \{0 \}$ is also positive, hence $0 \in B_i$, and from the hypothesis, there exists $q$ such that $B_q 0 \leqslant B_i$. This implies that for each $i$ such that $0 \in B_i$, as $A_i \times B_q 0 \leqslant p^* U$ one has $B_q A_i = p_!(A_i \times B_q 0) \leqslant U$ hence $A_i \triangleleft_q U$.

Now as $U \times \{0\}$ is locally positive and a subset of $p^*(U)$:

\[ U \times \{0\} \leqslant \bigvee_{i \in I \atop (A_i \times B_i) \wedge (U \times \{ 0 \})> \emptyset } \leqslant \bigvee_{i \in I \atop 0 \in B_i} A_i \times B_i \]

Applying $\pi_{1!}$ one gets (as any $B_i$ having a point is positive) that

\[U \leqslant \bigvee_{i \in I \atop 0 \in B_i} A_i \leqslant \bigvee_{i \in I \atop A_i \triangleleft U } A_i, \]

which concludes the proof of the first implication.

Assume now $(LB2)$, let $U$ be an arbitrary neighborhood of $0$, then as $\Hcal$ is metric, there exists an open sublocale $V$ such that $0 \in V$ and $V \triangleleft U$. In particular, there exists $q$ such that $B_q V \leqslant U$, and as $0 \in V$ one has:

\[ B_q 0 \subset B_q V \subset U \]
which proves $(LB1)$ and concludes the proof of the proposition.

}

}

\block{\Prop{Let $\Hcal$ be a pre-Banach locale, then its completion $\widetilde{\Hcal}$ is naturally endowed with a structure of Banach locale such that the map $\Hcal \rightarrow \widetilde{\Hcal}$ is a linear isometric map.}

\Dem{Everything comes more or less immediately from \ref{compextension} for the construction of operations and from \ref{strongdensity} and \ref{extensionofmetricprop} for the verification of the axioms:

Indeed, as $\Hcal \times \Hcal$ has a fiberwise dense image in $\widetilde{\Hcal} \times \widetilde{\Hcal}$, the canonical (uniform) map $p:\Hcal \times \Hcal \rightarrow \Hcal \rightarrow \widetilde{\Hcal}$ extends into a map $\widetilde{\Hcal} \times \widetilde{\Hcal} \rightarrow \widetilde{\Hcal}$. Similarly, the opposite map $m:\Hcal \rightarrow \Hcal$ is isometric and hence extends into a map $m: \widetilde{\Hcal} \rightarrow \widetilde{\Hcal}$ and one checks all the group axioms on $\widetilde{\Hcal}$ because they hold in $\Hcal$, that $\widetilde{\Hcal}$ is metric and that $\Hcal^n$ has a fiberwise dense image in $\widetilde{\Hcal}^n$.

The action of the locale of complex numbers on $\widetilde{\Hcal}$ is obtained in the same way: for each $\lambda \in Q[i]$ the multiplication by $\lambda$ is a uniform map $ \Hcal \rightarrow \Hcal $ and hence extends into a map $ \widetilde{\Hcal} \rightarrow \Hcal$, giving a map $\Q[i] \times \widetilde{\Hcal} \rightarrow \widetilde{\Hcal}$ and all the axioms of compatibility with the group law are also satisfied by a density argument.

Finally, we already know that there is a distance function on $\widetilde{\Hcal}$ we only have to check that $\Vert x \Vert = d(0,x)$ is a norm and that $d(x,y)=\Vert x - y \Vert$. But this also immediately comes from a density argument by \ref{extensionofmetricprop}.
}
}

\block{\Cor{Banach locale are complete metric locales.}
 
\Dem{Let $\Hcal$ be a Banach locale, in particular $\Hcal$ is a metric locale and hence by \ref{isometricsublocale} it identifies with a sublocale of $\widetilde{\Hcal}$. More precisely, as the inclusion is a linear map, $\Hcal$ identifies with a localic subgroup of a locally positive localic group $\widetilde{\Hcal}$, hence thanks to the constructive version of the closed subgroups theorem proved by P.T.Johnstone in \cite{johnstone1989constructive}, one concludes that $\Hcal$ is fiberwise closed (weakly closed in the terminology of \cite{johnstone1989constructive}) in $\widetilde{\Hcal}$ and hence is also complete (see the remark at the end of \ref{compdef}). }
}

\block{In particular, the action of $\Q[i]$ on a Banach locale extends to an action of its completion $\C$. Indeed (assuming that $\Hcal$ is complete), the map $B_n 0 \times \Q[i] \rightarrow \Hcal$ is uniform (it is $n$-Lipschitz) and hence it extends into $\overline{B_n 0} \times \C \rightarrow \Hcal$. One has a family of compatible maps $B_n 0 \times \C \rightarrow \Hcal$ which gives rise to a map $\Hcal \times \C \rightarrow \Hcal $. }

\block{Similarly to what is done in section \ref{Metricset}, a pre-Banach space in the usual (constructive) sense is exactly the same as a pre-Banach locale whose underlying locale is a discrete topological space. To such a Banach space one can associate its completion which is going to be a Banach locale. Conversely to any Banach locale one can associate its space of points which is a Banach space, and these two constructions induce an equivalence between the category of weakly spatial Banach locales (and linear map) and the category of Banach spaces (with bounded linear map).}

\subsection{The Localic Gelfand duality}

\block{\Def{A $C^*$ locale (or localic $C^*$ algebra) is a Banach locale $\Ccal$, endowed with an involution $* : \Ccal \rightarrow \Ccal$ and a product $\Ccal \times \Ccal \rightarrow \Ccal$ which satisfy the usual axioms for a $C^*$ algebra:

\begin{itemize}
\item $\Ccal$ is a $\C$ algebra (i.e. the product is associative, distributes over the addition and is compatible with the action of $\C$).
\item The $*$ involution is $\C$ anti-linear and satisfies $(ab)^*= b^* a^*$.
\item One has: $\Vert a b \Vert \leqslant \Vert a \Vert \Vert b \Vert$.
\item One has: $\Vert a^* a \Vert = \Vert a \Vert^2$.
\end{itemize}
}

All the axioms are equalities (or inequalities with respect to the specialization order), hence are clearly preserved by pull-back and therefore if $\Ccal$ is a $C^*$ algebra and $f$ is a geometric morphism to the base topos then $f^\#(\Ccal)$ is also a $C^*$ locale. And if $\Ccal$ is a (pre)-Banach locale endowed with an $*$ map and a map $\Ccal \times \Ccal \rightarrow \Ccal$ such that for some open surjection $f$, $f^\#(\Ccal)$ is a $C^*$ algebra for those structure then $\Ccal$ is a $C^*$ algebra.
}

\blockn{The main result of this section will be an anti-equivalence of categories between the categories of abelian unital $C^*$ locales and compact regular locales. The ``difficult" part lies in the construction of the two functors, and the proof that they are compatible with pull-back along geometric morphisms. Indeed once it is done, one can apply \ref{potentiallyweaklyspatial} to reduce the proof of the equivalence to the case of spatial $C^*$ algebras and completely regular compact locales which is already known (\cite{banaschewski2006globalisation} \cite{coquand2009constructive}). Actually, even the construction of the two functors could be avoided since we know that the notion of $C^*$ locale is the ``stackification" of the notion of $C^*$ algebra (it is a direct consequence of the observations made in \ref{stackification}), and one can prove (applying \ref{compregandwspa}) a similar result for compact regular locales and compact completely regular locales. Hence the already known equivalence between unital abelian $C^*$ algebras and compact completely regular locales immediately yields the equivalence between the ``stackified" notions, but we think that it is important to have an explicit construction of these functors without having to use descent theory.  }

\block{\Prop{Let $X$ be a compact regular locale, then $[X, \C]$ is a $C^*$ algebra, for the addition, product and involution given by the addition, the product and the complex conjugation of $\C$, and the norm given by:

\[ B_q0 = [X \ll f^* D_q] \]

where $D_q$ denotes the open disc of radius $q$ in $\C$, and $[X \ll f^* D_q]$ denotes the basic open which classifies the $f$ such that $X \ll f^* D_q$.

 }

\Dem{ $[X,\C]$ is indeed locally positive by \ref{X,Clocpos}. For the rest, we recall that Hyland gave in \cite{hyland1981function} a description of the theory classified by $[X,Y]$ in terms of the basic propositions $[U \ll f^* V]$ for $U \in \Ocal(X)$ and $V \in \Ocal(Y)$. From this description, we immediately obtain that:

\[ \bigvee_{q'<q} B_{q'} 0 = B_q 0;  \]
\[ \bigvee_{n} B_n 0 = [X,\C]. \]

Also, as $0$ is the point corresponding to the function constant equal to $0$, one has indeed $0 \in B_q 0$.

Hence the $B_q 0$ indeed define a function $\Vert . \Vert : [X,\C] \rightarrow \Length$ such that $\Vert 0 \Vert =0$, and such that $\bigvee_{n} B_n 0 = [X,\C]$.

All the algebraic axioms (including the triangular inequality) are checked on generalized point exactly as one does for classical points in the usual (constructive) case.

A basic open $[U \ll f^*V ]$ (for $U$ positive) contains $0$ if $U \ll \bigvee_{0 \in V} X$, but this implies that there exists a finite set $F$ included in $\{ 0 \in V \}$ such that $U \leqslant \bigvee_{f \in F} X$. A finite set is inhabited or empty, hence either $F$ is empty and $U = \emptyset$ or $F$ is inhabited and $0 \in V$. In the first case $[U \ll f^*V ] =[X,\C]$ contains all the $B_q 0$. In the second case one has a $q$ such that $D_q \ll V$ and hence $0 \in B_q 0= [X \ll f^*(D_q)] \leqslant [U \ll f^*(V)]$ which proves that the $B_q 0$ form a basis of neighborhood of $0$, and hence $[X,\C]$ is a Banach locale.

}
}

\block{\label{spectrum}We now want to construct the spectrum of a $C^*$ locale. We will start by defining the locale $\fn \Hcal$ of linear forms of norm smaller than $1$ on a Banach locale $\Hcal$ (the spectrum being the space of characters, it will be a sublocale of this locale). It generalizes the locale $\fn E$ constructed in \cite{mulvey1991globalization} and \cite{coquand1999direct}.

\Prop{Let $\Hcal$ be a Banach locale. There exists a sublocale $\fn \Hcal \subset [ \Hcal, \C ]_1$ which classifies the linear forms of norm smaller or equal to one on $\Hcal$. If $\Ccal$ is a unital commutative $C^*$ locale, then there exists a sublocale $\spec \Ccal \subset [\Ccal,\C]_1$ which classifies characters of $\Ccal$. }

\Dem{One can for example define the locale $\fn \Hcal$ as the intersection of the equalizer of the following two diagrams: 

\[ [ \Hcal, \C]_1 \rightrightarrows [  D_1 \times \Hcal, \C]_1 \]

where $D_1$ denotes the open unit ball in $\C$ and the two maps are the maps defined on generalized elements by: $f \mapsto \left ( (\lambda,x) \mapsto \lambda f(x) \right)$ and $f \mapsto \left( (\lambda,x) \mapsto f(\lambda x) \right)$, and where the distance on $D_1 \times \Hcal$ is the max distance.

And,

\[ [\Hcal, \C]_1 \rightrightarrows [ \Hcal \times \Hcal , \C]_1 \]

where $\Hcal \times \Hcal$ is endowed with the norm $\Vert x_1 \Vert+\Vert x_2 \Vert$ and the two maps are given by: $f \mapsto \left( (x,y) \mapsto f(x+y) \right)$ and $f\mapsto \left((x,y) \mapsto f(x)+f(y) \right)$.

\smallskip

A map $X \rightarrow \fn \Hcal$ is then exactly the data (internally to $X$) of a metric map from $\Hcal \rightarrow \C$ which is additive and linear with respect to complex numbers smaller than $1$. As it is also linear with respect to integers, it is linear on $n D_1$ for all $n$ and this forms an open cover of $\C$ so it concludes the proof.

\bigskip

If now $\Ccal$ is a unital $C^*$ locale, then one defines $\spec \Ccal$ as the intersection of the two previous equalizers with the pull-back of $\{1\} \subset \C$ by the map of evaluation on the unit on $[\Ccal, \C ]_1$ and with the equalizer of the following diagram:

\[ [\Ccal,\C]_1 \rightrightarrows [B_1 0 \times B_1 0, \C] \]

where $B_1 0$ is the open unit ball of $\Ccal$, and the distance $B_1 0 \times B_1 0$ is given by the max distance. The two maps are given by $f \mapsto \left( (x,y) \mapsto f(x)f(y) \right)$ and $f \mapsto \left( (x,y) \mapsto f(x y) \right)$.

A map factoring into $\spec \Ccal$ exactly corresponds to an internal character of $\Ccal$.

}

}

\block{The following  result is a localic version of the Banach-Alaoglu theorem.
\Prop{Let $\Hcal$ be a Banach locale, $\Ccal$ a unital commutative $C^*$ locale, then the locales $\fn \Hcal$ and $\spec \Ccal$ are compact regular locales.}

\Dem{ Compact regular locales descend along open surjections: for example because for a locale being compact and regular is the same thing as having a map to the point which is both proper and separated (see \cite{sketches} C.3.2.10) and because both proper maps and separated maps descend along open morphisms, (see \cite{sketches}C5.1.7). Hence it is enough to prove that some pull-back of $\fn \Hcal$ and $\spec \Ccal$ by an open surjection is compact and regular to conclude. In particular, by \ref{potentiallyweaklyspatial} one can freely assume that $\Hcal$ and $\Ccal$ are weakly spatial and hence that it is the completion of some Banach space $H$ or $C^*$ algebra $C$. But in this situation, a linear form or a character on the Banach locale is exactly the same as a linear form or a character on the set of points (by extension to the completion) and hence (the pull-back of) $\fn \Hcal$ and $\spec \Ccal$ classify the same theory as the locale $\fn H$ and $\spec C$ (also called $\text{MFn }C$) studied in \cite{mulvey1991globalization} and \cite{banaschewski2006globalisation} for the case of Grothendieck toposes, and in \cite{coquand1999direct} and \cite{coquand2009constructive} for general elementary toposes. These references prove that these locales are indeed compact (completely) regular.  }

 }

\block{\label{gelfand}\Th{The previous two constructions $X \rightarrow [X,\C]$ and $\Ccal \rightarrow \spec \Ccal$ induce an anti-equivalence of categories between unital abelian $C^*$ locales and compact regular locales.}

\Dem{These two constructions are defined in terms of the theory they classified and hence we can easily check that they are preserved by pull-back along geometric morphisms. They correspond to the well known notion of (completion of the) space of continuous functions on $X$ and spectrum of a $C^*$ algebra when $X$ is completely regular and when $\Ccal$ is weakly spatial. Moreover the two canonical maps ``evaluation at $x \in X$" from $X$ to $\spec [X,\C]$ and ``evalution at $c \in \Ccal$" from $\Ccal$ to $[\spec \Ccal,\C ]$ are preserved by pull-back (a proof by generalized points shows it immediately).

Hence, applying \ref{potentiallyweaklyspatial} one can pull-back (along an open surjection) those two maps to a similar situation but with $\Ccal$ and $[X,\C]$ weakly spatial (hence with $X$ completely regular by \ref{compregandwspa}). We can then conclude that the pull-back (along an open surjection) of the two canonical maps are isomorphisms from the usual constructive Gelfand duality (proved in \cite{banaschewski2006globalisation} for Grothendieck toposes, and generalized in \cite{coquand2009constructive} to arbitrary elementary toposes). And hence, as pull-back by an open surjection is conservative, these two canonical maps are isomorphisms. This proves that the two constructions are inverse from each other, the fact that they form an equivalence of categories follows immediately from the same argument.
}

}

\section*{References}

\bibliography{Biblio}

\begin{thebibliography}{10}

\bibitem{banaschewski2006globalisation}
{Banaschewski, Bernhard and Mulvey, Christopher J}.
\newblock {A globalisation of the Gelfand duality theorem}.
\newblock {\em {Annals of Pure and Applied Logic}}, 137(1):62--103, 2006.

\bibitem{borceux3}
{Borceux, F.}
\newblock {\em {Handbook of Categorical Algebra: Volume 3, Sheaf Theory}},
  volume~3.
\newblock {Cambridge University Press}, 1994.

\bibitem{bourbaki1966elements}
{Bourbaki, Nicolas}.
\newblock {\em {Elements of Mathematics: General Topology}}.
\newblock Hermann, 1966.

\bibitem{bunge1990descent}
{Bunge, Marta}.
\newblock An application of descent to a classification theorem for toposes.
\newblock In {\em {Mathematical Proceedings of the Cambridge Philosophical
  Society}}, volume 107, pages 59--79. {Cambridge Univ Press}, 1990.

\bibitem{burden1979banach}
{Burden, CW and Mulvey, CJ}.
\newblock {Banach spaces in categories of sheaves}.
\newblock In {\em Applications of sheaves}, pages 169--196. Springer, 1979.

\bibitem{coquand1999direct}
{Coquand, Thierry}.
\newblock {A Direct Proof of the Localic Hahn-Banach Theorem}.
\newblock {\em to appear}, 1999.

\bibitem{coquand2009constructive}
{Coquand, Thierry and Spitters, Bas and others}.
\newblock {Constructive Gelfand duality for C*-algebras}.
\newblock In {\em {Mathematical Proceedings of the Cambridge Philosophical
  Society}}, volume 147, pages 323--337. {Cambridge Univ Press}, 2009.

\bibitem{fell1977induced}
{Fell, James Michael Gardner and Douady, Adrien and Dal Soglio-H{\'e}rault,
  Letizia}.
\newblock {\em {Induced representations and Banach *-algebraic bundles}}.
\newblock Springer-Verlag Berlin, 1977.

\bibitem{mythesis}
Simon Henry.
\newblock {\em Des topos {\`a} la g{\'e}om{\'e}trie non commutative par
  l'{\'e}tude des espaces de {H}ilbert internes}.
\newblock PhD thesis, Universit{\'e} Paris 7, 2014.

\bibitem{hyland1981function}
{Hyland, JME}.
\newblock Function spaces in the category of locales.
\newblock In {\em Continuous lattices}, pages 264--281. Springer, 1981.

\bibitem{johnstone1989constructive}
{Johnstone, Peter T}.
\newblock A constructive “closed subgroup theorem” for localic groups and
  groupoids.
\newblock {\em {Cahiers de Topologie et G{\'e}om{\'e}trie Diff{\'e}rentielle
  Cat{\'e}goriques}}, 30(1):3--23, 1989.

\bibitem{sketches}
{Johnstone, P.T.}
\newblock {\em Sketches of an elephant: a topos theory compendium}.
\newblock Clarendon Press, 2002.

\bibitem{joyal1984extension}
{Joyal, A. and Tierney, M.}
\newblock {\em {An extension of the Galois theory of Grothendieck}}.
\newblock {American Mathematical Society}, 1984.

\bibitem{moerdijk1990classifying}
Ieke Moerdijk.
\newblock The classifying topos of a continuous groupoid. {II}.
\newblock {\em Cahiers de Topologie et G{\'e}om{\'e}trie Diff{\'e}rentielle
  Cat{\'e}goriques}, 31(2):137--168, 1990.

\bibitem{moerdijk1986connected}
{Moerdijk, Izak and Wraith, GC}.
\newblock Connected locally connected toposes are path-connected.
\newblock {\em {Transactions of the American Mathematical Society}},
  295(2):849--859, 1986.

\bibitem{mulvey1991globalization}
{Mulvey, Christopher J and Pelletier, Joan Wick}.
\newblock {A globalization of the Hahn-Banach theorem}.
\newblock {\em {Advances in Mathematics}}, 89(1):1--59, 1991.

\bibitem{picado2012frames}
{Picado, Jorge and Pultr, Ale{\'e}s}.
\newblock {\em {Frames and Locales: topology without points}}.
\newblock {Springer}, 2012.

\bibitem{vickers2005localic}
{Vickers, Steven}.
\newblock {Localic completion of generalized metric spaces I}.
\newblock {\em {Theory and Applications of Categories}}, 14(15):328--356, 2005.

\end{thebibliography}

\bibliographystyle{plain}

\end{document}